\documentclass[11pt]{article}
\usepackage{amssymb,amsmath,graphicx}
\usepackage{amsfonts,setspace,physics}
\usepackage{caption}
\usepackage[caption=false]{subfig}
\usepackage[numbers,sort,compress]{natbib}
\usepackage{epsfig,latexsym,graphicx,color}
\newtheorem{theorem}{Theorem}[section]

\newtheorem{proposition}[theorem]{Proposition}
\newtheorem{corollary}[theorem]{Corollary}
\newtheorem{remark}{Remark}[section]
\newtheorem{definition}{Definition}[section]

\begin{document}

\title{A Scale-invariant Generalization of the R\'{e}nyi Entropy,\\ 
Associated Divergences and their Optimizations\\
 under Tsallis' Nonextensive Framework}

\author{
Abhik Ghosh and Ayanendranath Basu\footnote{Corresponding author}\\
\small Interdisciplinary Statistical Research unit\\
\small Indian Statistical Institute, Kolkata, India  \\
\small{\it abhik.ghosh@isical.ac.in, ayanbasu@isical.ac.in}}
\maketitle

\begin{abstract}
Entropy and relative or cross entropy measures are two very fundamental concepts in information theory 
and are also widely used for statistical inference across disciplines. 
The related optimization problems, in particular the maximization of the entropy 
and the minimization of the cross entropy or relative entropy (divergence), are essential  for general logical inference in our physical world. 
In this paper, we discuss a two parameter generalization of 
the popular R\'{e}nyi entropy and associated optimization problems. 
We derive the desired entropic characteristics of the new generalized entropy measure 
including its positivity, expandability, {extensivity} and generalized (sub-)additivity. 
More importantly, when considered over the class of sub-probabilities, our new family turns out to be scale-invariant; 
this property does not hold for most existing generalized entropy measures. 
We also propose the corresponding cross entropy and relative entropy measures and discuss their geometric properties
including generalized Pythagorean results over $\beta$-convex sets.
The maximization of the new entropy and the minimization of the corresponding cross or relative entropy measures are 
carried out explicitly under the non-extensive (`third-choice') constraints given by the Tsallis' normalized $q$-expectations
which also correspond to the $\beta$-linear family of probability distributions.
Important properties of the associated forward and reverse projection rules are discussed along with their existence and uniqueness.
In this context, we have come up with, for the first time, a class of entropy measures 
-- a subfamily of our two-parameter generalization -- 
that leads to the classical (extensive) exponential family of MaxEnt distributions 
under the non-extensive constraints; this discovery has been illustrated through the useful concept of escort distributions
and can potentially be important for future research in information theory.
Other members of the new entropy family, however, lead to the power-law type generalized $q$-exponential MaxEnt distributions
which is in conformity with Tsallis' nonextensive theory. 
Therefore, our new family indeed provides a wide range of entropy and associated measures 
combining both the extensive and nonextensive MaxEnt theories under one umbrella.
\end{abstract}


\section{Introduction}
\label{SEC:concept}

The concept of entropy is a fundamental tool in information science, 
statistical physics and related statistical  applications.
Its development started at least one hundred years back in the context of thermodynamics and, interestingly,
this thermodynamic entropy increases over the `arrow of time' unlike all other physical variables 
making it a useful yet somewhat mysterious concept. 
Its wider application beyond thermodynamics, however,  started much later, 
after Shannon's groundbreaking  work leading to the development of mathematical information theory in 1948 \cite{Shannon:1948}
and Jaynes' presentation of the universal \textit{maximum entropy (MaxEnt) principle} 
for logical scientific inference in 1957 \cite{Jaynes:1957a, Jaynes:1957b}.
Shannon primarily defined the concept of information-theoretic entropy with 
the aim of developing an appropriate (uncertainty) measure of 
the amount of information lost in a noisy communication channel;
but a direct and highly interesting connection with the classical thermodynamic entropy 
can be made through Jaynes' MaxEnt principle (see \cite{Kapur/Kesavan:1992} for details).
Jaynes' work suggests that one should use all available information 
but be maximally uncommitted to the missing information 
leading to the (possibly constrained) maximization of Shannon's uncertainty measure; 
this initially provided a natural correspondence 
between statistical mechanics and general logical inference in information theory.
During the same period,  Kullback \cite{Kullback/Leibler:1951,Kullback:1952,Kullback:1953,Kullback:1954,Kullback:1956} 
provided the link between information theory  and Fisher's likelihood theory of general statistical inference 
along with the correspondence between a generalization of Jaynes' MaxEnt principle (minimum cross entropy principle) 
and Fisher's maximum likelihood principle; see also the monograph \cite{Kullback:1997} and the references therein.
These connections further enhanced the popularity of the entropy concept 
in a wide variety of fields of natural sciences \cite{Kapur:1990}.

In mathematical information theory, we consider the space of probability distributions over a finite alphabet set 
$\mathcal{A}=\{a_1, \ldots, a_n\}$ given by 
\begin{eqnarray}
\Omega_n = \left\{P=(p_1, \ldots, p_n) : p_i = \mbox{Prob}(a_i)\geq 0, \mbox{ for all }i=1, \ldots,n, W(P):=\sum_i p_i = 1 \right\},
\label{EQ:ProbSet_finite} 
\end{eqnarray}
and then define the Shannon entropy of any $P\in\Omega_n$  as
\begin{eqnarray}
\mathcal{E}^S(P) := -\sum_{i=1}^n p_i \ln p_i = \sum_{i=1}^n p_i \ln\left( \frac{1}{p_i}\right).
\label{EQ:Shannon_entropy}
\end{eqnarray}
Under certain desirable properties of the postulated uncertainty measures \cite{Shannon:1948},
Shannon obtained this unique form in (\ref{EQ:Shannon_entropy}), upto a positive multiplier, which he termed as the \textit{entropy}.
The quantity $I_i := \ln(1/p_i)$ is referred to as the \textit{elementary information gain}
associated with an event (alphabet) of probability $p_i$, or the \textit{code length} in information theory,
or the \textit{surprise} (less probable events are considered more ``surprising" than the more probable ones). 
Jaynes' MaxEnt principle suggests the prediction of the unknown natural distribution  
by maximizing $\mathcal{E}^S(P)$ over $P\in \Omega_n$ subject to the constraints of the given information,
which is used to model a wide range of systems (or situations) across disciplines. 

Subsequently, several generalizations  of the Shannon entropy have been developed for more complex systems;
among them the two most famous ones are the R\'{e}nyi entropy \cite{renyi:1961} and the Tsallis entropy \cite{Tsallis:1988}. 
Although the functional forms of these entropies were proposed much earlier,
their possible utility in explaining the behavior of natural systems 
was widely accepted only after more recent experimental verifications. 
In order to define these entropies in a more general context, 
let us consider the set of finite sub-probability distributions given by
\begin{eqnarray}
\Omega_n^\ast = \left\{P=(p_1, \ldots, p_n) : p_i \geq 0, \mbox{ for all }i=1, \ldots,n, W(P):=\sum_i p_i \leq 1 \right\} \supset \Omega_n.
\label{EQ:ProbSet_finite_gen} 
\end{eqnarray}
For any $P\in \Omega_n^\ast$, the corresponding (generalized) Shannon entropy turns out to be 
\begin{eqnarray}
\mathcal{E}^S(P) := -\frac{1}{W(P)}\sum_i p_i \ln p_i,
\label{EQ:Shannon_entropy_sp}
\end{eqnarray}
which coincides with  definition (\ref{EQ:Shannon_entropy}) for any $P\in\Omega_n$ (as $W(P)=1$);
see the derivation in \cite{renyi:1961}.
The R\'{e}nyi entropy of a general $P\in \Omega_n^\ast$ is defined, in terms of a tuning parameter $\alpha>0$, as
\begin{eqnarray}
\mathcal{E}_{\alpha}^R(P) := \frac{1}{1 - \alpha}\ln \left[\frac{\sum_i p_i^{\alpha}}{\sum_i p_i}\right],
~~~ \alpha> 0.
\label{EQ:renyi_entropy_sp}
\end{eqnarray}
For a probability distribution $P\in \Omega_n$, it further simplifies to 
\begin{eqnarray}\label{EQ:renyi_entropy}
\mathcal{E}_{\alpha}^R(P) &:=& \frac{1}{1 - \alpha}\ln \sum_i p_i^{\alpha} = \frac{\alpha}{1 - \alpha}\ln||P||_{\alpha},
~~~ \alpha>0,
\end{eqnarray}
where $||P||_\alpha$ denotes the $\alpha$-norm of  $P=(p_1, \ldots, p_n)$ defined as
$$
||P||_\alpha = \left(\sum\limits_{i=1}^n p_i^\alpha\right)^{\frac{1}{\alpha}}.
$$
In both (\ref{EQ:renyi_entropy_sp}) and (\ref{EQ:renyi_entropy}), 
the case $\alpha=1$ is defined through the limit as $\alpha\rightarrow 1$, 
which coincide with the Shannon entropy in (\ref{EQ:Shannon_entropy_sp}) and (\ref{EQ:Shannon_entropy}), respectively.
Interestingly, all members of this R\'{e}nyi entropy family are still extensive for the probability distributions, 
i.e., they satisfy 
$$
\mathcal{E}_{\alpha}^R(P\ast Q)=\mathcal{E}_{\alpha}^R(P)+\mathcal{E}_{\alpha}^R(Q), 
~~~\mbox{ for all } P\in \Omega_n, Q\in \Omega_m,~\alpha>0,
$$
where $P\ast Q=((p_iq_j))_{i=1, \ldots,n; j=1, \ldots, m}$ denotes the probability of 
the independent combination of two systems having probabilities $P$ and $Q$.

On the other hand, Tsallis entropy is the most popular non-extensive generalization of the Shannon entropy,
which is defined as
\begin{eqnarray}
\mathcal{E}_{q}^{T}(P) &:=& \frac{1-\sum_i p_i^{q}}{q-1} =   - \sum_i p_i^{q}\ln_q p_i,
~~~ q\in\mathbb{R}, ~P\in \Omega_n,
\label{EQ:Tsallis_entropy}
\end{eqnarray}
where $\ln_q$ denotes the deformed logarithm function (Definition \ref{DEF:q-deformed}).
The Tsallis entropy with index $q$ can be written as the $q$-deformed Shannon entropy
and coincides with the classical Shannon entropy (\ref{EQ:Shannon_entropy}) as $q\rightarrow 1$.
The quantity $q$ is also referred to as the nonextensivity index of the system, 
since we have the relation  
$$
\mathcal{E}_{q}^{T}(P\ast Q) = \mathcal{E}_{q}^{T}(P) + \mathcal{E}_{q}^{T}(Q) 
+ (1-q)\mathcal{E}_{q}^{T}(P)\mathcal{E}_{q}^{T}(Q),
~~~\mbox{ for all }  P\in \Omega_n,Q\in \Omega_m.
$$	
The theoretical prediction from this nonextensive entropy has later been 
seen to be extremely accurate for several advanced physical phenomena 
which leads to a whole new domain of nonextensive statistical framework \cite{Tsallis:2009}. 
One important component of the nonextensive framework is the generalization of the (linear) expectation constraints 
by the $q$-expectation or normalized $q$-expectation constraints \cite{Tsallis/etc:1998};
they are also related to the information theoretic concepts of escort distribution (Definition \ref{DEF:escort})
and $\beta$-linear family (Definition \ref{DEF:beta-linear})  which we will come back to in Section \ref{SEC:LNE_MaxEnt}.

Several further one and two-parameter generalizations of  the entropy functional have been proposed in the literature, 
although not all of them have equally significant applications with experimental validity.
We would like to mention two such generalizations of the R\'{e}nyi and Tsallis entropy, respectively,
known as Kapur's generalized entropy of order $\alpha$ and type $\beta$ \cite{Kapur:1967, Kapur:1969}
and the $(\alpha, \beta)$-norm entropy \cite{Joshi/Kumar:2016}. For any $P\in \Omega^\ast_n$, 
these families are defined, respectively, in terms of two positive (unequal) reals $\alpha$, $\beta$ as
\begin{eqnarray}
\mathcal{E}_{\alpha, \beta}^{K}(P) &:=& 
\frac{1}{\alpha-\beta}\ln\left(\frac{ \sum_i p_i^{\beta}}{\sum_i p_i^{\alpha}}\right),
\label{EQ:Kapur_entropy}\\
\mathcal{E}_{\alpha,\beta}^{N}(P) 
&:=& \frac{\alpha\beta}{\alpha-\beta}\left[\left({\sum_i p_i^\beta}\right)^{\frac{1}{\beta}}
- \left({\sum_i p_i^\alpha}\right)^{\frac{1}{\alpha}}\right] 
= \frac{\alpha\beta}{\alpha-\beta}\left[||P||_\beta - ||P||_\alpha\right].
\label{EQ:rsNorm_entropy}
\end{eqnarray}
Note that, for $P\in \Omega_n$, if we set either of the two tuning parameters $\alpha$ or $\beta$ to be one,
the first one reduces to the R\'{e}nyi entropy whereas the second one reduces to the Tsallis entropy. 
They both are symmetric with respect to $(\alpha, \beta)$ and are related in the limiting sense as
\begin{eqnarray}
\lim_{\alpha\rightarrow\beta}\mathcal{E}_{\alpha, \beta}^{K}(P)
= \frac{1}{\beta^2} \lim_{\alpha\rightarrow\beta}\mathcal{E}_{\alpha, \beta}^{N}(P)
= -\frac{\sum_i p_i^\beta \ln p_i}{\sum_i p_i^\beta}, ~~~P\in \Omega_n.
\label{EQ:AD_entropy}
\end{eqnarray}
The limiting functional in (\ref{EQ:AD_entropy}) is another a one-parameter family of generalized entropies
which had been previously studied independently by Aczel and  Daroczy \cite{Aczel/Daroczy:1963}
and will be referred to as the Aczel-Daroczy entropy $\mathcal{E}_{\beta}^{AD}(P)$;
note that $\mathcal{E}_{1}^{AD}(P)$ is the Shannon entropy in (\ref{EQ:Shannon_entropy}).

We emphasize the fact that neither the R\'{e}nyi nor  the Tsallis entropy
are scale-invariant over $P\in \Omega_n^\ast$; the same holds for their generalizations in (\ref{EQ:Kapur_entropy})--(\ref{EQ:AD_entropy}).
They are not even scale-equivariant except for the norm-entropy with $\alpha\neq \beta$.
This lack of invariance sometime makes the derivation of MaxEnt distribution and related statistics 
rather complicated while considering the general class of sub-probability distributions;
the MaxEnt distribution then exists only if $W(P)$ is pre-fixed (given).
Also, if we focus on measuring the pattern of the distribution only through the measure of uncertainty,
an appropriate entropy should be scale-invariant so that $P$ and $cP$ have the same entropy measure for any $c>0$
(as they have the same patterns over the state-space). 
In this paper, we will develop a new two-parameter family of entropy functionals,
generalizing the R\'{e}nyi entropy, that closely resembles the above two generalized families
but, in addition, provides the much desired scale-invariance property.  
We refer to our new generalized entropy as the logarithmic ($\alpha, \beta$)-norm entropy,
or simply the logarithmic norm-entropy (LNE); see Section \ref{SEC:LNE}.

Although this new LNE family has previously been introduced very briefly (as a generalized R\'{e}nyi entropy)
in our earlier work \cite{Ghosh/Basu:2018}  while describing a family of generalized relative entropy measures (and their applications),
its entropic characteristics and scope of application are practically unknown.
They will be developed here to justify its use as a general entropy functional.
Also, we will derive the maximum entropy (MaxEnt) theory for our LNEs under the non-extensive $q$-normalized expectation ($\beta$-linear family); 
the resulting MaxEnt distribution leads to a family of generalized exponential distributions 
(namely the $\beta$-power-law distributions in Definition \ref{DEF:beta-power-law}) having heavier tails. 
Our new MaxEnt theory resembles Tsallis' MaxEnt theory  yet generalizes it allowing a two-parameter structure with scale-invariance.

Another very important optimization problem related to entropy 
is the minimization of the associated cross entropy or relative entropy measures.
In information theory, we often have a prior guess of the distribution, say $Q$,
and the target distribution is then estimated through the minimization of a suitable cross entropy or a relative entropy measure from the prior $Q$. 
This is in direct correspondence with the MaxEnt principle;
in fact an appropriate minimizer of the cross entropy or the relative entropy (their forward projection; Definition \ref{DEF:Forward_Proj})
from a given uniform prior (i.e., no additional information) coincides with the MaxEnt distribution of the associated entropy measure. 
In statistics, the relative entropies are often referred to as  divergence measures and 
the minimization of appropriate divergences between the observed data and the postulated model  
(their reverse projection; Definition \ref{DEF:Reverse_Proj}) leads to robust inference  
in the presence of outliers or contamination in the data \cite{Basu/etc:2011}.
The most popular and classical relative entropy (divergence) measure related to the Shannon entropy is the 
Kullback-Leibler divergence (KLD) \cite{Kullback/Leibler:1951}
which is also related to maximum likelihood estimation in statistics \cite{Kullback:1997}.
For any two distributions $P,Q \in \Omega_n$, the KLD measure between them is given by
\begin{eqnarray}
\mathcal{RE}(P, Q) = \sum_i p_i\ln\left(\frac{p_i}{q_i}\right).
\label{EQ:RE_KL}
\end{eqnarray}
In the line of the R\'{e}nyi entropy, the R\'{e}nyi divergence generalizes the KLD measure; 
for  any two $P,Q \in \Omega_n$, their R\'{e}nyi divergence is defined as
\begin{eqnarray}
\mathcal{D}_\alpha(P, Q) = \frac{1}{\alpha-1} \ln\left(\sum_i p_i^\alpha q_i^{1-\alpha}\right),~~\alpha>0,~\alpha\neq 1.
\label{EQ:RE_renyi}
\end{eqnarray}
The R\'{e}nyi divergence for $\alpha=1$ is defined in the limiting sense and coincides with the KLD measure. 
However, unlike the KLD, the R\'{e}nyi divergences in (\ref{EQ:RE_renyi}) do not satisfy the important Pythagorean inequality 
over any given convex sets;
they only satisfy it over the $\alpha$-convex sets (Definition \ref{DEF:beta-convex}) as shown in \cite{VanErven/Harremos:2014}.
An alternative divergence measure, the relative $\alpha$-entropy, has been developed in \cite{Sundaresan:2007,Lutwak/etc:2005}
which satisfies the Pythagorean inequality over any given convex sets and is related to the R\'{e}nyi divergence 
via its definition. For any given $\alpha>0$ and $P,Q \in \Omega_n$, their relative $\alpha$-entropy can be defined as
\begin{eqnarray}
\mathcal{RE}_\alpha(P, Q) &:=&\mathcal{D}_{1/\alpha}(P_\alpha, Q_\alpha), 
\label{EQ:RE_alpha}
\end{eqnarray}
where $P_\alpha$ and $Q_\alpha$ are the $\alpha$-escort distributions (Definition \ref{DEF:escort}) of $P$ and $Q$, respectively.
Note that, at $\alpha=1$, the relative $\alpha$-entropy coincides with the R\'{e}nyi divergence of order 1 which is nothing but the KLD.
The geometric properties of this relative $\alpha$-entropy along with the corresponding minimization problems 
have been recently studied \citep{Kumar/Sundaresan:2015a,Kumar/Sundaresan:2015b}.

More recently, a two-parameter generalization of the relative $\alpha$-entropy measure has been studied in \cite{Ghosh/Basu:2018}
which is indeed related to our LNE measure and serves as the corresponding relative entropy (referred to as the LNRE). 
Additionally, in the present paper we will also develop a new cross entropy measure related to our LNE, to be referred to as the LNCE, 
and link it to the LNRE measure. 
The  geometric properties of the LNCE and LNRE measures will be discussed along with the associated Pythagorean-type results. 
The minimization problems (projection rules) associated with these generalized information measures (LNCE and LNRE)
will also be studied under the Tsallis' non-extensive constraints (equivalently under the $\beta$-linear family).

In brief, we summarize the main contributions of this manuscript as follows:
\begin{itemize}
	\item We propose and study the detailed properties of a new class of entropy measures 
	which is scale-invariant over the space of (finite) sub-probability distributions
	and contains the popular R\'{e}nyi entropy class (and hence also the Shannon entropy) on the 
	space of  probability distributions (having finite support). 
	In particular, we prove that the new LNE family satisfies all the entropic characteristic axioms of 
	R\'{e}nyi entropy except possibly the generalized additivity property (Definition \ref{DEF:Gen_additivity}). 
	Other than the R\'{e}nyi subfamily, only one subfamily of LNE at $\alpha=\beta$ satisfies the generalized additivity property; 
	all other members of the LNE family are shown to have a generalized sub-additivity property with suitably chosen weight-functions.
	However, like other non-Shannon entropies, LNE does not satisfy  the branching or the recursivity property (Definition \ref{DEF:Branching}).
	
	\item Along with its scale-invariant nature, our LNE family is also shown to be extensive in nature; 
	for any two independent (communication or physical) systems over a finite alphabet set, 
	the LNE value (entropy) of their combination equals the sum of their individual LNE measures.
	Like the R\'{e}nyi entropy, this important extensivity property of our LNE measures 
	makes them useful in the context of information theory to analyze complex extensive systems.

	\item We derive the MaxEnt distribution corresponding to the new LNE family 
	under the Tsallis non-extensive constraint (of the third kind) given in terms of the normalized $q$-expectation,
	or equivalently under the $\beta$-linear families. 
	When the two parameters of the LNE family differ, the resulting MaxEnt distributions 
	form the generalized exponential family of the $\beta$-power-law distributions (Definition \ref{DEF:beta-power-law}) having heavier tails;
	this resembles and generalizes the MaxEnt theory and applications corresponding to the Tsallis and the R\'{e}nyi entropies.
	More interestingly, the new LNE subfamily at $\alpha=\beta$ provides the usual exponential family of MaxEnt 
	distribution even under the non-extensive constraint. 
	Indeed, this is the \textit{first} family of generalized entropies in the literature (as per the knowledge of the authors) 
	which is shown to provide the usual Shannon-type MaxEnt theory under the Tsallis non-extensive framework.
	
	\item The LNE family can be linked with (and can also be motivated from) 
	the generalized relative ($\alpha, \beta$)-entropy measure (LNRE) studied in \cite{Ghosh/Basu:2018}.
	In the present paper, we further define the corresponding   cross entropy measure (LNCE) connecting them
	and completing the full circle of important measures in information theory.
	Several important geometric properties of these new families of LNCE and LNRE measures are noted.
	In the particular case $\beta=1$, the LNCE between two probability measures
	coincides with the R\'{e}nyi divergence in (\ref{EQ:RE_renyi})
	whereas the corresponding LNRE measures coincide with the relative $\alpha$-entropy in (\ref{EQ:RE_alpha}).
	A generalized Pythagorean property is proved for these new LNRE and LNCE measures over the $\beta$-convex sets,
	which extends the corresponding results for R\'{e}nyi entropy presented in \cite{VanErven/Harremos:2014,Kumar/Sason:2016}, 
	indicating their further utility in information theory.

	\item Projection problems are one of the basic tools in mathematical information theory 
	when a prior distributional guess is known for the given information problem. 
	In this paper, we also define and study the projection rules arising from the proposed LNRE and LNCE.
	The forward projection rules are seen to be the same for both LNCE and LNRE measures and
	a set of sufficient conditions for their existence and uniqueness are discussed for general cases.
	Due to the lack of invariance of LNCE in its second argument, the reverse projection is discussed only for the LNRE measures.

	\item As a particular example, the forward projection of the LNRE (and LNCE) measures over the $\beta$-linear family is explored in great detail.
	An explicit form of the resulting projection is derived and their properties are studied;
	we have verified if and when this projection rule satisfies the  important characteristic axioms of \cite{Csiszar:1991}
	over the $\beta$-convex set of probability distributions.  
	This forward projection problem indeed also corresponds to the minimization of the LNRE (or LNCE) under 
	the Tsallis non-extensive constraints of the  third kind.
	Our derivations extend the corresponding results for the minimum KLD distribution in \cite{Dukkipati/etc:2005} 
	obtained under non-extensive or extensive frameworks, respectively, for the LNREs (or LNCEs) with $\alpha\neq \beta$ or $\alpha=\beta$. 
	Again, we get an interesting subfamily of relative entropy and cross entropy measures at $\alpha=\beta$,
	which leads to extensive (Shannon-type) results under the non-extensive frameworks	and 
	provides a potentially new direction  of research combining the two concepts.

	\item The LNE can also be interpreted as a suitable R\'{e}nyi entropy 
	of the escort distribution associated with any given sub-probability distribution.
	Similar relations are also derived for our LNRE (and LNCE) with the  R\'{e}nyi divergence measure.
	Noting that the escort distributions generate a one-to-one correspondence \cite{Karthik/Sundaresan:2018}, 
	the basic information geometric properties of the new LNE and LNRE measures are seen to be equivalent 
	with the corresponding results for the existing R\'{e}nyi entropy  and divergence measures. 
	In particular, the Pythagorean property and the forward projection rules generated by the LNRE (or LNCE)
	are shown to be equivalent with those generated by the R\'{e}nyi divergence. 
	Such one-to-one correspondences clearly justify the usefulness of the proposed LNE and associated LNRE/LNCE
	in extending the R\'{e}nyi information concepts  with the additional feature of scale-invariance,
	potentially leading to better results in applications to complex systems.
	
	\item As a particular application, the usefulness of the LNRE in robust statistical inference has been demonstrated in detail, 
	with appropriate discussions, mathematical justifications, statistical insights, illustrations  and references. 
	Given a set of independent and identically distributed data, 
	the parameter of a assumed statistical model family can be obtained by the reverse LNRE projection of the data distribution 
	(relative frequencies in case of discrete data) onto the assumed parametric family. 
	The resulting minimum LNRE estimators (MLNREEs), with appropriate choice of the two defining tuning parameters, 
	are shown to provide better trade-offs between the (asymptotic) efficiency  under pure data
	and robustness under data contamination compared to other related (existing) minimum divergence estimators.  
\end{itemize}

For brevity in presentation and quick reference, 
all the important existing concepts, referred to previously and also throughout the rest of the paper,
are collectively defined in an Appendix to the paper.

\section{A Two-parameter Generalization of the R\'{e}nyi Entropy}
\label{SEC:LNE}

We consider first the set $\Omega_n^\ast$ of all sub-probability distributions over the finite state-space 
as defined in (\ref{EQ:ProbSet_finite_gen}). Given two positive reals $\alpha, \beta$, 
we define a generalized entropy measure of any $P\in \Omega_n^\ast$ as
\begin{eqnarray}
\mathcal{E}_{\alpha, \beta}^{LN}(P) 
:= \frac{\alpha\beta}{\alpha-\beta}\ln\left[\frac{\left({\sum_i p_i^\beta}\right)^{\frac{1}{\beta}}}{
	\left({\sum_i p_i^\alpha}\right)^{\frac{1}{\alpha}}}\right] 
= \frac{\alpha\beta}{\alpha-\beta}\left[\ln||P||_\beta - \ln||P||_\alpha\right] ,
~~~ \alpha\neq\beta.
\label{EQ:LN_entropy}
\end{eqnarray}
We can extend its definition at $\alpha=\beta$ through the limiting functional as $\alpha\rightarrow\beta$, 
which is given by 
\begin{eqnarray}
\mathcal{E}_{\beta, \beta}^{LN}(P) 
:= -\beta\frac{\sum_i p_i^\beta \ln p_i}{\sum_i p_i^\beta}  + \ln\left(\sum_i p_i^\beta\right)
= \beta\left[\mathcal{E}_{\beta}^{AD}(P) + \ln||P||_\beta\right],
~~~ \beta>0.
\label{EQ:LN_entropy0}
\end{eqnarray}
Note that $\mathcal{E}_{1,1}^{LN}(P)$ does not necessarily coincide with the corresponding (generalized) Shannon entropy
(\ref{EQ:Shannon_entropy_sp}) for a sub-probability $P\in \Omega_n^\ast$, 
but does so for the probability distributions having unit sum ($P\in \Omega_n$).
Similarly, if $P\in \Omega_n$ and any one of the two tuning parameters $\alpha$, $\beta$ assumes the value one, 
the generalized entropy in (\ref{EQ:LN_entropy}) coincides with the R\'{e}nyi entropy in (\ref{EQ:renyi_entropy}).


The major advantage of the functional forms of the our proposed entropy, given in (\ref{EQ:LN_entropy}) and (\ref{EQ:LN_entropy0}),
is its scale-invariance property: $\mathcal{E}_{\alpha, \beta}^{LN}(cP) =\mathcal{E}_{\alpha, \beta}^{LN}(P)$
for any $P\in \Omega_n^\ast$ and $c, \alpha, \beta >0$ such that $cP\in\Omega_n^\ast$. 
This striking property is satisfied neither by the Shannon entropy nor its existing generalizations  
like the R\'{e}nyi entropy, the Tsallis entropy or those given in (\ref{EQ:Kapur_entropy})--(\ref{EQ:rsNorm_entropy}).
However, there are clear similarities of the proposed functional form in (\ref{EQ:LN_entropy})
with the entropy formulas in (\ref{EQ:Kapur_entropy})--(\ref{EQ:rsNorm_entropy}) which motivated our construction. 
Instead of considering the difference between two norms as in the $(\alpha, \beta)$-Norm entropy in (\ref{EQ:rsNorm_entropy}),
we consider the difference between their logarithms which generates a scale-invariant functional form. 
Similarly the difference with the entropy in  (\ref{EQ:Kapur_entropy}) can be seen only in the power of the terms within the logarithm
which is imposed to make our entropy in (\ref{EQ:LN_entropy}) scale-invariant.  

Therefore, it appears that the proposed entropies in (\ref{EQ:LN_entropy})--(\ref{EQ:LN_entropy0})  are 
the first two-parameter generalizations of  the Shannon and R\'{e}nyi entropy over $P\in\Omega_n$
that is scale-invariant over the larger set of sub-probability distributions $\Omega_n^\ast$.
The functional forms (\ref{EQ:LN_entropy})--(\ref{EQ:LN_entropy0}) were initially found in \cite{Ghosh/Basu:2018}
where the authors had referred to them as the possible generalized R\'{e}nyi entropies (without any detailed study). 
However, since there are several generalized versions of the R\'{e}nyi entropy available in the literature,
by noting their similarity with the $(\alpha, \beta)$-Norm entropy in (\ref{EQ:rsNorm_entropy}), 
we will denote this generalized family of entropies given by (\ref{EQ:LN_entropy}) and (\ref{EQ:LN_entropy0})
as the \textit{Logarithmic ($\alpha, \beta$)-Norm entropy}, or simply the \textit{Logarithmic Norm-Entropy (LNE)}, of $P\in\Omega_n^\ast$.
It is interesting to note that $\mathcal{E}_{\alpha, \beta}^{LN}(P)$ is symmetric in the tuning parameters $\alpha, \beta$.

\subsection{Examples and Limiting Cases}

Before undertaking a detailed theoretical investigation of the LNE family, 
let us start with some simple yet important examples to demonstrate its behavior. 

\begin{figure}[!h]
	\centering
	\subfloat[$\alpha=0.1$]{
		\includegraphics[width=0.3\textwidth]{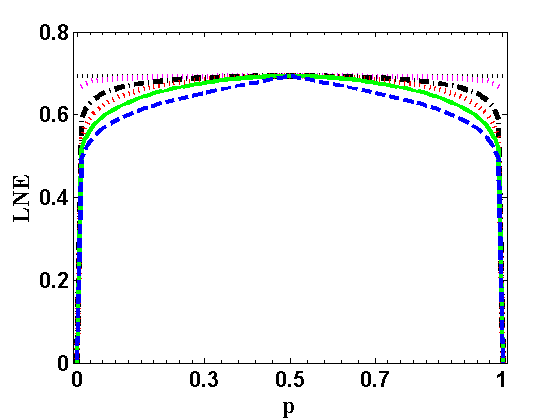}
		\label{FIG:9size_unknown_100_0}}
	~ 
	\subfloat[$\alpha=0.5$]{
		\includegraphics[width=0.3\textwidth]{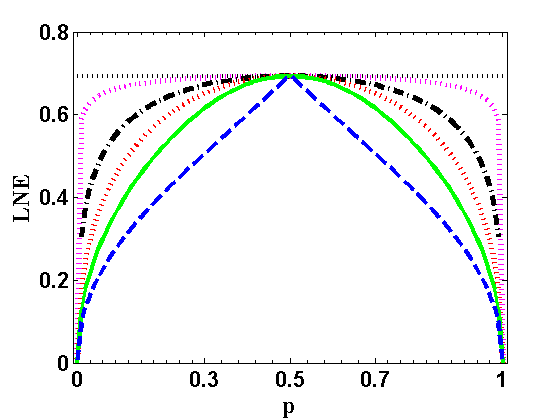}
		\label{FIG:9size_unknown_100_05}}
	~ 
	\subfloat[$\alpha=1$]{
		\includegraphics[width=0.3\textwidth]{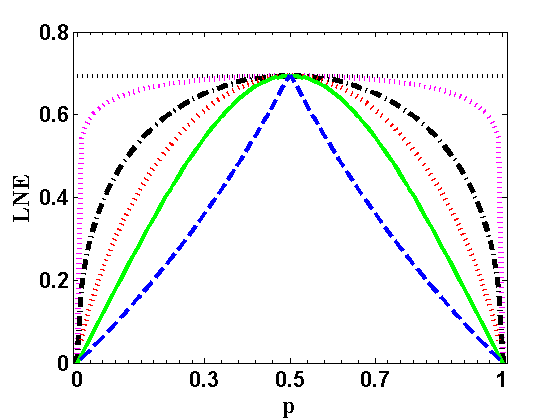}
		\label{FIG:9size_unknown_100_1}}
	\\
	\subfloat[$\alpha=2$]{
		\includegraphics[width=0.3\textwidth]{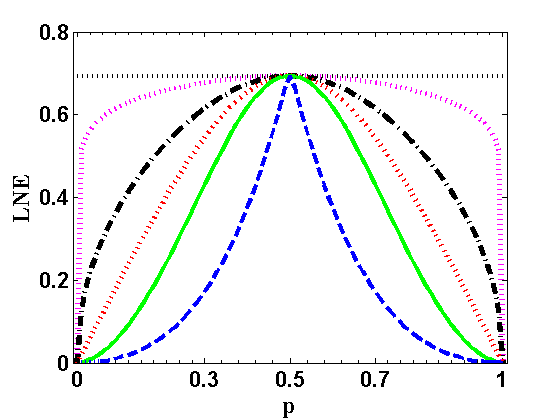}
		\label{FIG:9power_unknown_100_0}}
	~ 
	\subfloat[$\alpha=10$]{
		\includegraphics[width=0.3\textwidth]{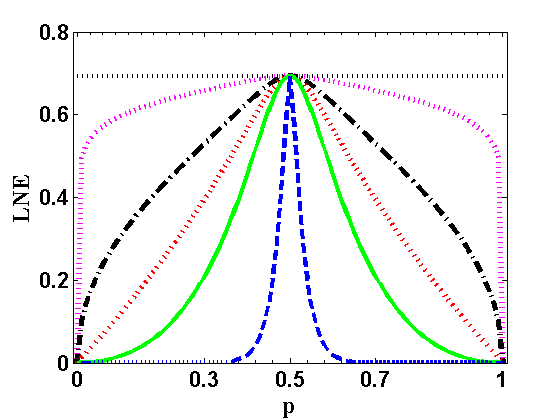}
		\label{FIG:9power_unknown_100_05}}
	~ 
	\subfloat[$\alpha=100$]{
		\includegraphics[width=0.3\textwidth]{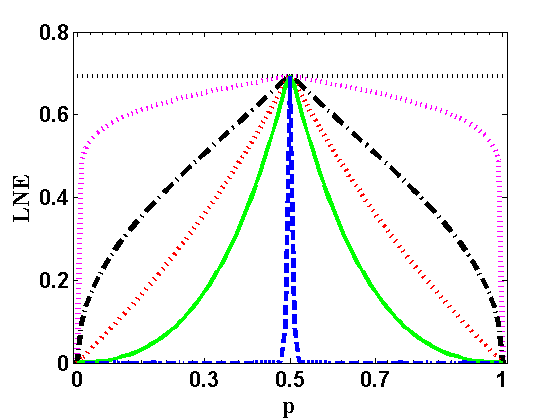}
		\label{FIG:9power_unknown_100_1}}
	\caption{Values of the LNE of $P=\{p, 1-p\}$ plotted against $p\in[0,1]$ for different values of $\alpha, \beta$.
		\footnotesize[$\beta=0$: black dotted; $\beta=0.1$: magenta dotted; $\beta=0.5$: black dash-dotted; 
		$\beta=1$: red dotted; $\beta=2$: green solid; $\beta=100$: blue dashed lines]	}
	\label{FIG:LNE_Ber}
\end{figure}

\bigskip
\noindent
\textbf{Example 1. }[\textit{LNE of the Bernoulli Distribution}]\\
Consider the Bernoulli distribution which has two states ($n=2$)
and is characterized by the success rate $p$ of any one state.
In Figure \ref{FIG:LNE_Ber}, we have plotted the LNE $\mathcal{E}_{\alpha, \beta}^{LN}(\{p, 1-p\})$
over the success probability $p\in [0,1]$  for  different values of $\alpha, \beta$.
Clearly, as expected, 
all members of the LNE family attain their  minimum (zero) and maximum ($\ln(2)=0.693$) at $p=0,1$ (degenerate distributions) 
and $p=1/2$ (uniform distribution), respectively. 
Also they are all continuous in $p\in[0,1]$. 

Additionally, we can observe that the LNE values decrease as any one of $\alpha$ or $\beta$ increases
from zero while the other is kept fixed. 
The forms of the LNE for such limiting cases are derived in the following theorem for general (sub-)probability distributions.
\hfill{$\square$}
\\

We now derive the forms of the LNE for some limiting cases which is presented in the following theorem. 

\begin{theorem}
	\label{THM:LNE_prop3} 
	Take any $P\in \Omega_n^\ast$ and $\beta>0$. 
	\begin{enumerate}
		\item[a)] As $\alpha\rightarrow 0$, $\mathcal{E}_{\alpha, \beta}^{LN}(P)\rightarrow \ln(n)$, 
		the maximum entropy value, independently of $\beta$ and $P$.
		\item[b)] As  $\alpha\rightarrow \infty$, 
		$\mathcal{E}_{\alpha, \beta}^{LN}(P)\rightarrow \beta\left[-\ln(p_{\max})+\ln||P||_\beta\right]$,
		which can be thought of as a scale-invariant generalization of the Min-entropy given by $-\ln(p_{\max})$.
		Here $p_{\max} = \max_i p_i$.	
	\end{enumerate}
\end{theorem}  
\noindent\textbf{Proof:}
For a given $P\in \Omega_n^\ast$ and $\beta>0$, we can rewrite the LNE at $\alpha\neq \beta$ as 
\begin{eqnarray}
\mathcal{E}_{\alpha, \beta}^{LN}(P) 
:= \frac{\alpha}{\alpha-\beta}\ln\left({\sum_i p_i^\beta}\right) + \frac{\beta}{\beta-\alpha}\ln\left({\sum_i p_i^\alpha}\right).
\label{EQ:LNE_limit}
\end{eqnarray}
(a) Taking limit as $\alpha\rightarrow 0$, the first term in (\ref{EQ:LNE_limit}) converges to zero, 
whereas the second term converges to $\ln(n)$.
\\
(b) Taking limit as $\alpha\rightarrow\infty$, the first term $\frac{1}{1-(\beta/\alpha)}\ln\left(\sum_i p_i^\beta\right)$ in (\ref{EQ:LNE_limit})
tends to  $\ln\left(\sum_i p_i^\beta\right)$. But, the second term in (\ref{EQ:LNE_limit}) is of the form ($\frac{-\infty}{-\infty}$)
as $\alpha\rightarrow\infty$, and hence we use L'Hospital rule to get 
\begin{eqnarray}
\lim\limits_{\alpha\rightarrow\infty} \frac{\beta}{\beta-\alpha}\ln\left({\sum_i p_i^\alpha}\right)
= \lim\limits_{\alpha\rightarrow\infty} \frac{\beta}{-1}\frac{\sum_i p_i^\alpha\ln p_i}{\sum_i p_i^\alpha}
= - \beta \ln(p_{\max}).
\end{eqnarray} 
Combining the limits of both the terms, we get the desired result.
\hfill{$\square$}

\bigskip
The above theorem indicates the nature of the LNE over its tuning parameters when one of them is fixed finitely. 
Note that the scale-invariant generalization obtained at $\alpha\rightarrow \infty$ represents 
the Min-entropy of the escort measure of $P$. 
We conjecture, based on our empirical examinations, that the value of $\mathcal{E}_{\alpha, \beta}^{LN}(P)$ 
monotonically decreases as $\alpha$ increases from 0 to $\infty$, 
at least for most common probability distributions $P$ if not for all of $\Omega_n^\ast$.

Next, to get an idea about their behaviors when both $\alpha, \beta$ vary simultaneously,
let us study the LNEs of the binomial distribution.

\begin{figure}[!b]
	\centering
	\subfloat[$n=10, p=0.1$]{
		\includegraphics[width=0.3\textwidth]{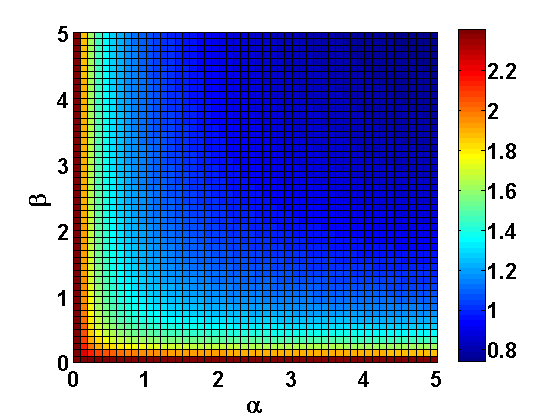}
		\label{FIG:9size_unknown_100_0}}
	~ 
	\subfloat[$n=10, p=0.3$]{
		\includegraphics[width=0.3\textwidth]{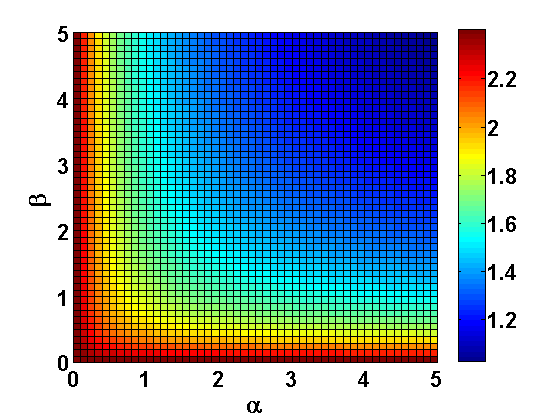}
		\label{FIG:9size_unknown_100_05}}
	~ 
	\subfloat[$n=10, p=.5$]{
		\includegraphics[width=0.3\textwidth]{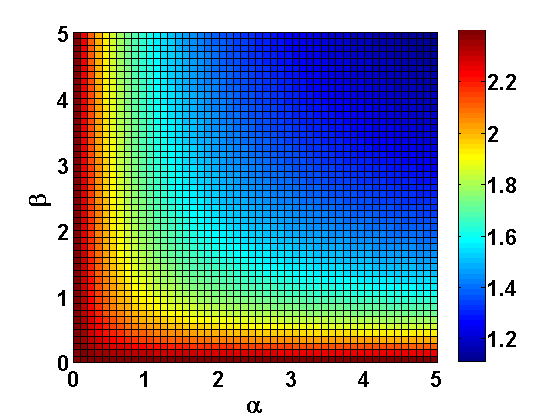}
		\label{FIG:9size_unknown_100_1}}
	\\
	\subfloat[$n=100, p=0.1$]{
		\includegraphics[width=0.3\textwidth]{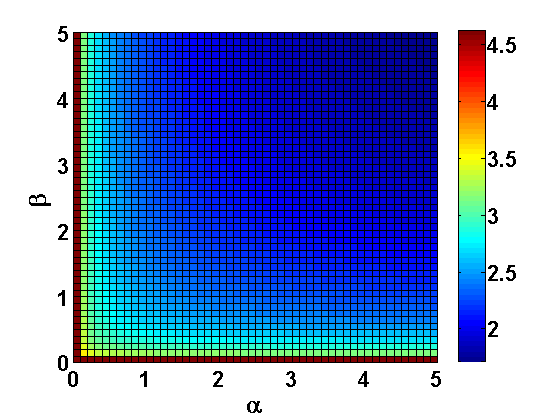}
		\label{FIG:9power_unknown_100_0}}
	~ 
	\subfloat[$n=100, p=0.5$]{
		\includegraphics[width=0.3\textwidth]{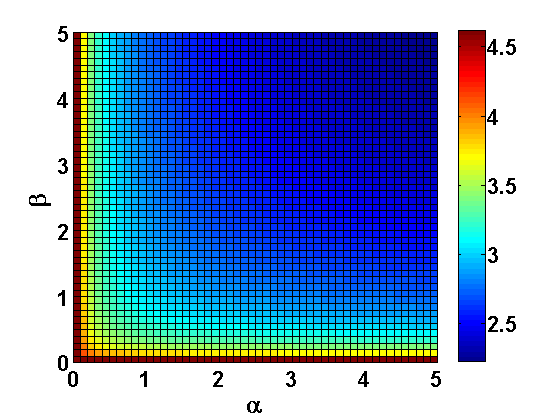}
		\label{FIG:9power_unknown_100_05}}
	~ 
	\subfloat[$n=100, p=0.9$]{
		\includegraphics[width=0.3\textwidth]{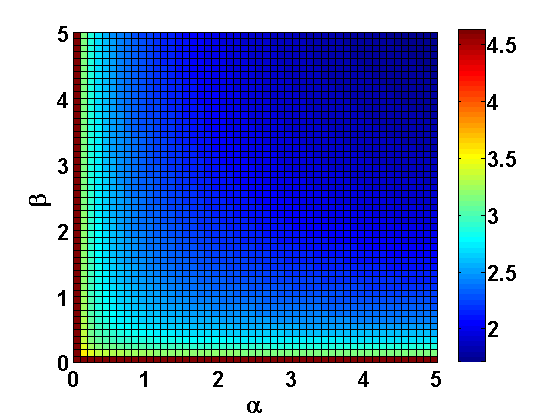}
		\label{FIG:9power_unknown_100_1}}
	\caption{Vales of the LNE of different Bin($n,p$) distribution plotted against $(\alpha, \beta)$.}
	\label{FIG:LNE_Bin}
\end{figure}

\bigskip\noindent
\textbf{Example 2. }[\textit{LNEs of the Binomial Distribution}]\\
Consider $n$ states having the binomial probability structure with success rate $p$;
this situation arises quite frequently in information theory while communicating an $n$-bit information over a noisy channel. 
We have computed and plotted, in Figure \ref{FIG:LNE_Bin}, 
the LNE values of this binomial distribution over $(\alpha, \beta)$ for different $n$ and $p$.
We can see that, for any binomial distribution, the entropy is maximized at $(\alpha, \beta)\rightarrow (0, 0)$
and decreases further as $(\alpha, \beta)$ moves away from zero (towards positive infinity).
Also, for a fixed $\alpha, \beta$ and a fixed number ($n$) of the state-space,
the LNE is maximized over the family of binomial distributions at $p=1/2$ 
and symmetrically decreases in either side leading to the minimum value of zero at $p=0$ and $p=1$. 
\hfill{$\square$}
\\

\subsection{Entropic Characteristics}

The Shannon as well as the R\'{e}nyi entropy are obtained via some desired axiomatic characteristics 
that a entropy must satisfy as a measure of uncertainty. 
Any new entropy measure must also satisfy those axiomatic postulates or their appropriate generalizations
which we now verify for our LNE.
The first result justifies that the LNE can indeed be considered as a measure of uncertainty. 

\begin{proposition}[\cite{Ghosh/Basu:2018}]
For any $\alpha, \beta>0$, the LNE functionals defined in (\ref{EQ:LN_entropy})--(\ref{EQ:LN_entropy0}) 
are always non-negative for all $P\in \Omega_n^\ast$.
They equal zero for the degenerate distributions and take the maximum value $\ln(n)$, over $\Omega_n^\ast$, 
if and only if all $p_i$s are equal (i.e., $P$ is uniform).
\label{THM:LNE_prop1}
\end{proposition}

Interestingly, unlike other generalizations of the R\'{e}nyi entropy, 
the maximum value of an LNE which is attained at the uniform distribution 
is independent of the tuning parameters $(\alpha, \beta)$ and is the same as that of the classical Shannon entropy.
Hence the LNE family provides a universal framework of comparison with a fixed bounded range of entropy values, 
namely $[0, \ln(n)]$, providing, at the same time, different structures to explain different types of systems 
through two tuning parameters $\alpha, \beta>0$.
In addition, the maximum value of the entropy increases further as the number of (microscopic) states ($n$) increases, as desired. 

The next theorem verifies if the members of this LNE family satisfy some additional desired 
characteristics of the R\'{e}nyi or other generalized entropies over its domain $\Omega_n^\ast$ or $\Omega_n$.

\begin{theorem}
\label{THM:LNE_prop2} 
For any $P\in \Omega_n^\ast$ and $\alpha, \beta>0$, the LNE $\mathcal{E}_{\alpha, \beta}^{LN}(P)$ 
satisfies the following properties:
\begin{enumerate}
	\item[a)] $\mathcal{E}_{\alpha, \beta}^{LN}(P)$ is continuous in the $p_i$s for all $i$ with $p_i\geq 0$ (unless all $p_i=0$).
	\item[b)] $\mathcal{E}_{\alpha, \beta}^{LN}(P)$ is a symmetric function of $(p_1, \ldots, p_n)$.
	\item[c)] $\mathcal{E}_{\alpha, \beta}^{LN}(\{1,0\})=\mathcal{E}_{\alpha, \beta}^{LN}(\{0,1\})=0$.
	~~{\normalfont [Decisivity]}
	\item[d)] If $P=\{p\}\in \Omega_1^\ast$ for any $p\in(0,1]$, then $\mathcal{E}_{\alpha, \beta}^{LN}(\{p\})=0$.
	\item[e)] For any $P=(p_1, \ldots, p_n)\in \Omega_n^\ast$, we have 
$\mathcal{E}_{\alpha, \beta}^{LN}(P)=\mathcal{E}_{\alpha, \beta}^{LN}(\{p_1, \ldots, p_n, 0\})$.
~~{\normalfont [Expandability]}\\
Note that, here $(p_1, \ldots, p_n, 0)$ is an element of $\Omega_{n+1}^*$.
	\item[f)] For $P=(p_1, \ldots, p_n)\in \Omega_n^\ast$ and $Q=(q_1, \ldots, q_m)\in\Omega_m^\ast$, 
	let us define their independent  combination as $P\ast Q = (p_iq_j)_{i=1, \ldots, n; j=1, \ldots, m}$.
	Then, 
	$$
	\mathcal{E}_{\alpha, \beta}^{LN}(P\ast Q)=\mathcal{E}_{\alpha, \beta}^{LN}(P)+\mathcal{E}_{\alpha, \beta}^{LN}(Q).
	~~~\mbox{{\normalfont [Shannon additivity/Extensivity]}}
	$$
	\item[g)] $\mathcal{E}_{\alpha, \beta}^{LN}(P)$, at any $n\geq 2$,  does not satisfy the branching 
or the recursivity properties (unlike the Shannon entropy) as defined in Definition \ref{DEF:Branching}.
\end{enumerate}
\end{theorem}  
\noindent\textbf{Proof:}\\
(a) The proof for the case $\alpha\neq \beta$ follows directly from the continuity of the norm functionals $||P||_\alpha$, $||P||_\beta$ 
and the logarithmic function, 
whereas the proof of the $\alpha=\beta$ case follows from the continuity of the Aczel-Daroczy entropy \cite{Aczel/Daroczy:1963} 
and the norm functional $||P||_\beta$.\\
(b--c) These two properties follow directly from the definition of LNE.\\
(d) Note that, by the definition of the norm, $||\{p\}||_{\gamma} = p$ for all $\gamma> 0$ and $p\in(0,1]$. 
Hence, for $\alpha\neq \beta$, we get 
$$
\mathcal{E}_{\alpha, \beta}^{LN}(\{p\}) = \frac{\alpha\beta}{\alpha-\beta}[\ln(p) - \ln(p)] = 0.
$$
Also, for $\alpha=\beta$, we have 
$$
\mathcal{E}_{\beta, \beta}^{LN}(\{p\}) = \beta\left[- \frac{p^\beta \ln(p)}{p^\beta} + \ln(p)\right] = 0.
$$

\noindent
(e) It follows trivially from definitions, since $||\{p_1, \ldots, p_n, 0\}||_{\gamma} = ||P||_{\gamma}$ for all $\gamma> 0$
and the Aczel-Daroczy entropy is extendable \cite{Aczel/Daroczy:1963}. \\
(f) First note that, for any $\gamma>0$, we have
$$
||P\ast Q||_{\gamma} = \left(\sum_{i=1}^n\sum_{j=1}^m (p_iq_j)^\gamma\right)^{1/\gamma}
= \left(\sum_{i=1}^np_i^\gamma \sum_{j=1}^m q_j^\gamma\right)^{1/\gamma} = ||P||_\gamma \cdot ||Q||_{\gamma}.
$$
Therefore, for $\alpha\neq \beta$ ($\alpha, \beta>0$), we get 
\begin{eqnarray}
\mathcal{E}_{\alpha, \beta}^{LN}(P\ast Q) &=& \frac{\alpha\beta}{\alpha-\beta}[\ln||P\ast Q||_\beta - \ln||P\ast Q||_\alpha] 
\nonumber\\
&=& \frac{\alpha\beta}{\alpha-\beta}[\ln||P||_\beta + \ln||Q||_\beta  - \ln||P||_\alpha - \ln||Q||_\alpha]
=  \mathcal{E}_{\alpha, \beta}^{LN}(P) + \mathcal{E}_{\alpha, \beta}^{LN}(Q).
\nonumber
\end{eqnarray}
For $\alpha=\beta$, on the other hand, we can use the Shannon additivity of the Aczel-Daroczy entropy \cite{Aczel/Daroczy:1963} to get 
\begin{eqnarray}
\mathcal{E}_{\beta, \beta}^{LN}(P\ast Q) &=& \beta[\mathcal{E}_{\beta}^{AD}(P\ast Q) + \ln||P\ast Q||_\alpha] 
\nonumber\\
&=& \frac{\alpha\beta}{\alpha-\beta}[\mathcal{E}_{\beta}^{AD}(P)  + \mathcal{E}_{\beta}^{AD}(Q)  + \ln||P||_\alpha + \ln||Q||_\alpha]
=  \mathcal{E}_{\beta, \beta}^{LN}(P) + \mathcal{E}_{\beta, \beta}^{LN}(Q).
\nonumber
\end{eqnarray}
(g) Finally, to show that the LNE does not satisfy the branching or the recursivity properties, let us consider the simple probabilities 
$P=\{1-p, p\}$ and $Q=\{q, 1-q\}$ for some $p,q\in[0,1]$.
Clearly, at $\alpha\neq \beta$, we have
\begin{eqnarray}
&&\mathcal{E}_{\alpha, \beta}^{LN}(\{1-p, pq, p(1-q)\})
\nonumber\\
&=&\frac{\alpha\beta}{\alpha-\beta}\left[\frac{1}{\beta}\ln\left((1-p)^\beta + p^\beta (q^\beta + (1-q)^\beta)\right) 
- \frac{1}{\alpha}\ln\left((1-p)^\alpha + p^\beta (q^\alpha + (1-q)^\alpha)\right) \right]
\nonumber\\
&\neq& \frac{\alpha\beta}{\alpha-\beta}\left[\ln||P||_\beta + p^a \ln||Q||_{\beta}
- \ln||P||_\alpha - p^a \ln||Q||_{\alpha}\right]
= \mathcal{E}_{\alpha, \beta}^{LN}(P) + p^a \mathcal{E}_{\alpha, \beta}^{LN}(Q), 
\nonumber
\end{eqnarray}
for any $a>0$. Similarly, the case $\alpha=\beta$ can be proved which is skipped for brevity.
\hfill{$\square$}

\bigskip
As the R\'{e}nyi entropy is characterized by the generalized-mean \cite{Kolmogorov:1930, Nagumo:1930} additivity property 
(Definition \ref{DEF:Gen_additivity}), it is important to check the same for our proposed extension as well. 
However, the only subclasses of the LNE family satisfying this property 
are the R\'{e}nyi entropy class and the family in (\ref{EQ:LN_entropy0}), although with different weight functions.
Other members of the LNE family are sub-additive in the same generalized mean with an 
appropriately chosen weight function, as shown in the following theorem.

\begin{theorem}[Generalized-Mean Sub-additivity]
For any two sub-probability distributions $P=(p_1, \ldots, p_n)\in \Omega_n^\ast$ and $Q=(q_1, \ldots, q_m)\in\Omega_m^\ast$
with $W(P) + W(Q) \leq 1$,  let us define the combined system (sub)-probability $P\cup Q = (p_1, \ldots, p_n, q_1, \ldots, q_m)$
and take $g(x)=2^{\left(1-\frac{\alpha}{\beta}\right)\frac{x}{c}}$ for $c=\ln2$ and pre-fixed $\alpha,\beta>0$ with $\alpha\neq \beta$,
$\alpha\neq 1$, $\beta\neq 1$. 
Then, we have
\begin{eqnarray}
\mathcal{E}_{\alpha, \beta}^{LN}(P\cup Q) 
&\leq& 	g^{-1}\left[\frac{||P||_\beta^\alpha g\left(\mathcal{E}_{\alpha, \beta}^{LN}(P)\right)
	+||Q||_\beta^\alpha g\left(\mathcal{E}_{\alpha, \beta}^{LN}(Q)\right)}{
	||P||_\beta^\alpha+||Q||_\beta^\alpha}\right], ~~~~~~~~\mbox{if } ~0<\alpha<\beta.
\label{EQ:GN_mean_Subadditivity}
\end{eqnarray}
For $\alpha > \beta>0$, the inequality in (\ref{EQ:GN_mean_Subadditivity}) is reversed. \\
\label{THM:LNE_propGM}
\end{theorem}  
\noindent\textbf{Proof: }
Fix $P\in \Omega_n^\ast$, $Q\in \Omega_m^\ast$ and $\beta>0$ 
as described in the statement of the theorem and take any $\alpha>0$.
By definition of the LNE family for $\alpha\neq \beta$, one can deduce
$$
g(\mathcal{E}_{\alpha, \beta}^{LN}(P)) = \frac{||P||_\alpha^\alpha}{||P||_\beta^\alpha} 
= \frac{\sum_i p_i^\alpha}{\left(\sum_i p_i^\beta \right)^{\frac{\alpha}{\beta}}},
$$
and similarly for $Q$. 
Now, the condition $W(P) + W(Q) \leq 1$ ensures that  $P\cup Q\in \Omega_n^\ast$ and hence we get
\begin{eqnarray}
g(\mathcal{E}_{\alpha, \beta}^{LN}(P\cup Q)) 
&=& \frac{||P||_\alpha^\alpha}{||P||_\beta^\alpha} 
= \frac{\sum_i p_i^\alpha + \sum_i q_i^\alpha}{\left(\sum_i p_i^\beta + \sum_i q_i^\beta\right)^{\frac{\alpha}{\beta}}}
=\frac{||P||_\beta^\alpha g(\mathcal{E}_{\alpha, \beta}^{LN}(P)) + ||Q||_\beta^\alpha g(\mathcal{E}_{\alpha, \beta}^{LN}(Q)) }{
	\left(||P||_\beta^\beta+||Q||_\beta^\beta\right)^{\frac{\alpha}{\beta}}}.
\nonumber
\end{eqnarray}
Further, by the order property of $L_p$-norm of the vector $(||P||_\beta, ||Q||_\beta)$, we get
$$
\left(||P||_\beta^\beta+||Q||_\beta^\beta\right)^{\frac{1}{\beta}} 
\geq 
\left(||P||_\beta^\alpha+||Q||_\beta^\alpha\right)^{\frac{1}{\alpha}},
~~~~~~\mbox{ if }~0<\alpha<\beta,
$$  
where the inequality is reversed for $0<\beta<\alpha$ and becomes equality at $\alpha=\beta$.
Combining the above two relations, we get the sub-additivity result (\ref{EQ:GN_mean_Subadditivity}) at $\alpha\neq \beta$. 
\hfill{$\square$}

\bigskip
Note that the right-hand side of Equation (\ref{EQ:GN_mean_Subadditivity}) represents a generalized-mean 
of the LNE values of $P$ and $Q$ defined through the link function $g$ and weights $(||P||_\beta^\alpha, ||Q||_\beta^\alpha)$.
By the symmetry of the LNE family with respect to the two tuning parameters $\alpha, \beta$,
we can always obtain the general sub-additivity property (\ref{EQ:GN_mean_Subadditivity}),
extending from  Definition \ref{DEF:Gen_additivity}, for any member of the LNE family with a suitable choice of weights; 
the respective weights will be $(||P||_\beta^\alpha, ||Q||_\beta^\alpha)$ or $(||P||_\alpha^\beta, ||Q||_\alpha^\beta)$,
according to $\alpha>\beta$ or $\alpha<\beta$.
Also, if we define the LNE in terms of logarithm base 2, as in the case of R\'{e}nyi entropy, 
we have $c=\log_2(2)=1$ in the link function $g$ in Theorem \ref{THM:LNE_propGM}.
For the cases $\alpha=1$ or $\beta=1$, our LNE coincides with the R\'{e}nyi entropy and hence satisfy the generalized additivity property 
as per Definition \ref{DEF:Gen_additivity}.


\subsection{Correspondence with the R\'{e}nyi Entropy via Escort Distributions}
\label{SEC:renyi_LNE}

We have already noted that the R\'{e}nyi entropy is a special class of the proposed LNE family
and they both have many similar entropic characteristics. 
Another interesting interpretation of the new LNE family can be observed 
through the so-called \textit{escort distribution} \cite{Beck/Schogl:1995,Abe:2003,Chapeau-Blondeau/etc:2011}
that puts it in one-to-one correspondence with the class of R\'{e}nyi entropies 
and provides further insights and justifications for our LNE.

The $\beta$-escort distribution for any $P\in \Omega_n^\ast$ and given $\beta>0$ is formally defined in Definition \ref{DEF:escort}.
For any subset $\mathbb{E}\subseteq \Omega_n^\ast$, 
let us denote the set of corresponding $\beta$-escort distributions by 
$$
\mathbb{E}^{(\beta)} = \left\{ P_\beta : P\in \mathbb{E} \right\} \subseteq \Omega_n.
$$
It can easily be shown that \cite[e.g.,][]{Karthik/Sundaresan:2018} for any $\mathbb{E}\subseteq \Omega_n^\ast$, 
with $W(P)$ being constant for every $P\in \mathbb{E}$, there is a one-to-one correspondence between the sets $\mathbb{E}$ and $\mathbb{E}^{(\beta)}$
via the $P \leftrightarrow P_\beta$ correspondence (Definition \ref{DEF:p-pbeta}).
In particular, the set of probability distributions $\Omega_n$ is in one-to-one correspondence with 
the set of their escort distributions $\Omega_n^{(\beta)}$.   
The following proposition, obtained based on this correspondence, would be extremely useful in the subsequent part of the paper;
the proof follows in the same line of Propositions 1 and 2 of \cite{Karthik/Sundaresan:2018} and is hence omitted.

\begin{proposition}
For any subset $\mathbb{E}\subseteq \Omega_n$ and any $\beta>0$, we have the following results.
\begin{itemize}
	\item[a)] A functional is concave (or convex) over the set $\mathbb{E}^{(\beta)}$ if 
	and only if it is $\beta$-concave (or $\beta$-convex) over $\mathbb{E}$ as per Definition \ref{DEF:beta-concaveF}. 
	\item[b)] The set $\mathbb{E}^{(\beta)}$ is convex  if and only if $\mathbb{E}$ is $\beta$-convex   as per Definition \ref{DEF:beta-convex}. 
	\item[c)] The set $\mathbb{E}$ is convex if and only if $\mathbb{E}^{(\beta)}$ is $(1/\beta)$-convex.
	\item[d)] The set $\mathbb{E}$ is closed in $L_\beta$ norm if and only if $\mathbb{E}^{(\beta)}$ is closed in $L_1$ norm. 
\end{itemize}
\label{PROP:escort}
\end{proposition}

Now, for any $P\in\Omega_n$, we can re-express the LNE given in (\ref{EQ:LN_entropy}) and (\ref{EQ:LN_entropy0}) alternatively as 
\begin{eqnarray}
\mathcal{E}_{\alpha, \beta}^{LN}(P) 
= \mathcal{E}_{\alpha/\beta}^R(P_{\beta}) = \mathcal{E}_{\beta/\alpha}^R(P_{\alpha}),
~~~~
\mathcal{E}_{\beta, \beta}^{LN}(P) = \mathcal{E}^S(P_{\beta}).
\label{EQ:LNE_escort}
\end{eqnarray}
Therefore, the newly proposed LNE is nothing but the R\'{e}nyi entropy of order $(\alpha/\beta)$ 
for the corresponding $\beta$-escort distribution. 
Since there is a one-to-one relation between $(\alpha, \beta)$ and $(\gamma=\alpha/\beta, \beta)$, 
the above relation also provides a one-to-one correspondence between the LNE values of the probability distributions in $\Omega_n$
with the R\'{e}nyi entropy values of the (escort) probabilities in  $\Omega_n^{(\beta)}$.
This equivalence can further be extended to the LNE of sub-probability distributions within $\Omega_n^\ast$
by utilizing the scale-invariance property of LNE. For any $P\in\Omega_n^\ast$, its LNE value 
$\mathcal{E}_{\alpha, \beta}^{LN}(P)$ equals $\mathcal{E}_{\alpha, \beta}^{LN}(P') $ with $P'= P/W(P) \in\Omega_n$. 
Therefore, for any $\alpha, \beta>0$,  our proposed LNE values over $\Omega_n^\ast$
are in direct one-to-one correspondence with the R\'{e}nyi entropy values of the class of escort distributions.
More generally, we can say that 
$$
\left\{\mathcal{E}_{\alpha, \beta}^{LN}(P) : P\in\Omega_n^\ast, \alpha>0, \beta>0 \right\} 
= \left\{ \mathcal{E}_{\gamma}^R(P) : P\in\Omega_n^{(\beta)}, \gamma>0, \beta>0  \right\}.
$$

The above correspondences clearly justify the use of the LNE functionals as a generalized class of entropy functionals.
Besides the equivalence with the values of R\'{e}nyi entropies, our proposed LNEs provide a much larger class of entropy functions
that can be used to model different varieties of systems via its two tuning parameters.
Therefore, the LNE family has the potential to provide more flexibility in all the applications where the R\'{e}nyi entropy is traditionally used;
further discussions are provided later in Sections  \ref{SEC:renyi_world}-\ref{SEC:Applications}.

\subsection{Concavity}

An important property of an entropy measure is its concavity that facilitates the maximum entropy theory.
It is fortunate for usual entropies like Shannon, R\'{e}nyi and Tsallis that they turn out to be concave 
which allows us to restrict ourselves only to the search of  a local maximum (under any constraint);
any local maximum will be a global one by their concavity. The following theorem examines this 
important property of our proposed LNE family for suitably chosen values of the tuning parameters $(\alpha, \beta)$.

\begin{theorem}
\label{THM:LNE_prop4} 
Suppose that either of the following two conditions on $(\alpha, \beta)$ holds.
\begin{itemize}
	\item[a)] $0<\beta \leq 1$ and $\alpha\geq \beta$ is such that $\ln||P||_\alpha$ is convex in $P$.
	\item[b)] $0<\alpha \leq 1$ and $\beta\geq \alpha$ is such that $\ln||P||_\beta$ is convex in $P$.
\end{itemize} 
Then, the LNE $\mathcal{E}_{\alpha, \beta}^{LN}(P)$ is concave in $P\in \Omega_n^\ast$.
\end{theorem}  
\noindent\textbf{Proof:}
We will prove the theorem under Condition (a). Then, it also holds under Condition (b) by symmetry of LNE in $(\alpha, \beta)$.\\
So, let us assume the Condition (a) holds and take $P, Q\in \Omega_n^\ast$, $\lambda\in [0,1]$. 
Since $\beta\leq 1$, by Minkowski inequality, we have
$$
||\lambda P + (1-\lambda) Q||_\beta \geq \lambda ||P||_\beta + (1-\lambda)||Q||_\beta.
$$ 
Combining it with the monotonicity and concavity of logarithmic function, we get 
$$
\ln||\lambda P + (1-\lambda) Q||_\beta \geq \ln\left[\lambda ||P||_\beta + (1-\lambda)||Q||_\beta\right] 
\geq \lambda \ln||P||_\beta + (1-\lambda)\ln||Q||_\beta.
$$ 
On the other hand, by convexity of $\ln||P||_\alpha$, we get 
$$
\ln||\lambda P + (1-\lambda) Q||_\alpha \leq \lambda \ln||P||_\alpha + (1-\lambda)\ln||Q||_\alpha.
$$ 
Thus, along with $\alpha>\beta$, we finally get 
\begin{eqnarray}
\mathcal{E}_{\alpha, \beta}^{LN}(\lambda P + (1-\lambda)Q) 
&=& \frac{\alpha\beta}{\alpha - \beta} \left[\ln||\lambda P + (1-\lambda) Q||_\beta - \ln||\lambda P + (1-\lambda) Q||_\alpha \right]
\nonumber\\
&\geq &  \frac{\alpha\beta}{\alpha - \beta} \left[  \lambda \ln||P||_\beta + (1-\lambda)\ln||Q||_\beta 
- \lambda \ln||P||_\alpha - (1-\lambda)\ln||Q||_\alpha\right]~~
\nonumber\\
&=& \frac{\alpha\beta}{\alpha - \beta} \left[  \lambda \left\{\ln||P||_\beta - \ln||P||_\alpha\right\} 
+ (1-\lambda)\left\{\ln||Q||_\beta - \ln||Q||_\alpha\right\}\right]
\nonumber\\
&=& \lambda \mathcal{E}_{\alpha, \beta}^{LN}(P) + (1-\lambda) \mathcal{E}_{\alpha, \beta}^{LN}(Q).  
\end{eqnarray}
This proves the concavity of LNE under Condition (a).
\hfill{$\square$}

\begin{remark}
Note that the concavity of the LNEs requires additional condition on the values of the tuning parameters $(\alpha, \beta)$. 
This assumption is clearly satisfied for the case $\alpha=1$ or $\beta=1$ (the R\'{e}nyi divergence). 
Additionally it is trivially satisfied as either of $\alpha$ or $\beta$ tends to zero or infinity.
However, it still remains an open problem to verify this assumption for more general and arbitrary collection of distributions
for other values of $\alpha$ and $\beta$. 
\end{remark}

Whenever the assumption of Theorem \ref{THM:LNE_prop4} is not satisfied, the LNEs would still be concave over 
the set of escort distributions, via their correspondence with the R\'{e}nyi entropy in (\ref{EQ:LNE_escort}), for all $\alpha, \beta >0$.
This result, presented in the following theorem, will suffice to study the MaxEnt distributions for LNEs
via the $P \leftrightarrow P_\beta$ correspondence.

\begin{theorem}[Weaker Concavity Results for LNE]
\label{THM:LNE_prop5} ~~
\begin{itemize}
	\item[a)] If $0<\alpha\leq \beta$, the LNE functional $\mathcal{E}_{\alpha, \beta}^{LN}(\cdot)$ is concave over 
	the set of escort distributions $\Omega_{n}^{(\beta)}$ 	and hence it is $\beta$-concave over the set of probability distributions $\Omega_{n}$. 
	\item[b)] If $0< \beta\leq \alpha$, the LNE functional $\mathcal{E}_{\alpha, \beta}^{LN}(\cdot)$ is concave over 
	the set of escort distributions $\Omega_{n}^{(\alpha)}$
	and hence it is $\alpha$-concave over the set of probability distributions $\Omega_{n}$. 
	\item[c)] For any $\alpha, \beta>0$, the LNE functional $\mathcal{E}_{\alpha, \beta}^{LN}(\cdot)$ is pseudoconcave and Schur concave over 
	either set of escort distributions $\Omega_{n}^{(\beta)}$ or $\Omega_{n}^{(\alpha)}$.
\end{itemize} 
\end{theorem}  
\noindent\textbf{Proof:}
The first results in Part (a) follows from the correspondence (\ref{EQ:LNE_escort}) and the concavity of R\'{e}nyi entropy of order $\leq 1$ 
\cite{Ben-Bassat/Raviv:1978}. 
The second result follows from Part (a) of Proposition \ref{PROP:escort} invoking scale-invariance of the LNEs.

\noindent Results in Part (b) follow from those of (a) by using the fact that the LNE does not change 
while interchanging the tuning parameters $\alpha, \beta$.

\noindent Results in part (c) can again be obtained using the correspondence (\ref{EQ:LNE_escort}) 
along with the properties of the R\'{e}nyi entropy  from \cite{Ben-Bassat/Raviv:1978,Ho/Verdu:2015}. 
\hfill{$\square$}

\section{The Maximum Entropy Distribution under Tsallis' Nonextensive Constraint}
\label{SEC:LNE_MaxEnt}

The maximum entropy principle is a fundamental concept in inferential science including information theory, statistics and statistical physics.
We will now derive the maximum entropy (MaxEnt) distribution corresponding to the new LNE family under appropriate sets of constraints.
%
It has been observed, more recently,  that 
the  generalized entropies like the Tsallis or R\'{e}nyi entropy provide more accurate predictions 
under the non-extensive constraints given in terms of the normalized $q$-expectation
instead of the linear expectation \cite{Tsallis/etc:1998, Kumar/Sason:2016}. 
So, in order to obtain the MaxEnt distribution,  here also we consider 
a set of $m$ non-extensive constraints given by
\begin{eqnarray}
\frac{\sum_i g_r(a_i)p_i^q}{\sum_i p_i^q} = G_r, ~~~r=1, \ldots, m,
\label{EQ:Constraints_Nq}
\end{eqnarray}
where $g_1, \ldots, g_m$ are given functions on $\mathcal{A}$ and $G_1, \ldots, G_m$ are fixed constants.
Note that, if $P=(p_1, \ldots, p_n)\in \Omega_n$ 
and $q=1$, the left hand side of the constraints in (\ref{EQ:Constraints_Nq}) simplifies to the usual (extensive) linear expectation;
otherwise it known as the normalized $q$-expectation of $g_r$ (denoted as $\expval{\expval{g_r}}_q$).
As noted formally in \cite{Kumar/Sason:2016}, the set of probability distributions satisfying the  constraints in (\ref{EQ:Constraints_Nq})
with $q=\beta$ corresponds to the $\beta$-linear family (Definition \ref{DEF:beta-linear}), 
a notion often used in modern generalized information theory \cite{Karthik/Sundaresan:2018,Kumar/Sason:2016b,Gayen/Kumar:2019}.

For the maximization of the LNE having parameters $\alpha, \beta$, by symmetry,
we can consider $q$ to be either of these two parameters. 
For concreteness, in this paper, we will consider $q=\beta$
and $(\alpha, \beta)$ satisfies the conditions of Theorem \ref{THM:LNE_prop4} or \ref{THM:LNE_prop5}
so that a maximizer of the LNE functional exists.
Then, 
by the scale-invariance property, it is enough to maximize  the LNE functional over 
the larger set $\Omega_n^\ast$ and then the required MaxEnt distribution over $\Omega_n$ can be obtained 
simply by normalizing the maximizer functional over $\Omega_n^\ast$.
This makes the maximization process practically easier and 
we can either consider the arguments of standard calculus or those of functional analysis and variational calculus. 
The final MaxEnt distribution obtained under the constraints in (\ref{EQ:Constraints_Nq}) 
is derived in the following theorem for all $\alpha, \beta>0$. 

\begin{theorem}
The probability distribution $P=(p_1, \ldots, p_n)\in\Omega_n$ that 
maximizes the LNE functional $\mathcal{E}_{\alpha, \beta}^{LN}(\cdot)$  
under the non-extensive constraints in (\ref{EQ:Constraints_Nq}) with $q=\beta$ 
(i.e., maximizes it over a $\beta$-linear family generated by the functions $(g_r - G_r)$, $r=1, \ldots, m$)
has the following forms. 
\begin{itemize}
	\item[a)]  When $\alpha\neq \beta$, there exists constants $\lambda_1, \ldots, \lambda_m$ such that 
	\begin{eqnarray}
	{p}_i &=& \frac{1}{Z} \left[1 + \left(\alpha - {\beta}\right)\sum_{r=1}^m \lambda_r (g_r(a_i)-G_r)\right]^{\frac{1}{\alpha-\beta}},
	~~~~~i=1, \ldots, n,
	\label{EQ:MaxEnt_final}\\
	Z &=& \sum_{i=1}^n \left[1 + \left(\alpha - {\beta}\right)\sum_{r=1}^m \lambda_r (g_r(a_i)-G_r)\right]^{\frac{1}{\alpha-\beta}}.
	\end{eqnarray}
	
	\item[b)] When $\alpha=\beta$, there exists constants $\lambda_1, \ldots, \lambda_m$ such that
	\begin{eqnarray}
	{p}_i &=& \frac{\exp\left[\sum_{r=1}^m \lambda_r (g_r(a_i)-G_r)\right]}{
		\sum_i \exp\left[\sum_{r=1}^m \lambda_r (g_r(a_i)-G_r)\right]},
	~~~~~i=1, \ldots, n.
	\label{EQ:MaxEnt_final_0}
	\end{eqnarray}
In both cases, the $m$ constants  $\lambda_1, \ldots, \lambda_m$ are uniquely defined by the $m$ constraints in (\ref{EQ:Constraints_Nq}). 
\end{itemize} 
\label{THM:LNE_MaxEnt}
\end{theorem}
\noindent\textbf{Proof:}
Based on the discussions preceding the theorem, to further simplify the optimization problem, let us define $\widetilde{P} = P/||P||_\beta$
for any $P\in \Omega_n^\ast$. Then, from its definition and (\ref{EQ:Constraints_Nq}), we have 
\begin{eqnarray}
\sum_i \widetilde{p}_i^\beta=1, ~~~~~
\sum_i g_r(a_i)\widetilde{p}_i^\beta=G_r, ~~r=1, \ldots,m. 
\label{EQ:Constraint_simple}
\end{eqnarray}
Since $\mathcal{E}_{\alpha, \beta}^{LN}(P)=\mathcal{E}_{\alpha, \beta}^{LN}(\widetilde{P})$
by scale-invariance, it is enough to maximize $\mathcal{E}_{\alpha, \beta}^{LN}(\widetilde{P})$
subject to the restriction given by (\ref{EQ:Constraint_simple}).
Using the methods of Lagrange multiplier, our objective function is now given by 
\begin{eqnarray}
F(\widetilde{P}) = \mathcal{E}_{\alpha, \beta}^{LN}(\widetilde{P}) +\lambda_0\left(\sum_i \widetilde{p}_i^\beta-1\right)
+\sum_{r=1}^m \lambda_r \left(\sum_i g_r(a_i)\widetilde{p}_i^\beta-G_r\right),
\label{EQ:Obj_Fun}
\end{eqnarray}
where $\lambda_r$s are Lagrange multipliers for $r=0, 1, \ldots, m$.

\noindent\textbf{Part (a):} \\
Let us first assume $\alpha\neq\beta$. Then, 
the first order condition for optimization of $F(\widetilde{P})$ is given by
\begin{eqnarray}
\frac{\alpha\beta}{\alpha-\beta}\left[\frac{\widetilde{p}_i^{\beta-1}}{\sum_i \widetilde{p}_i^{\beta}}
- \frac{\widetilde{p}_i^{\alpha-1}}{\sum_i \widetilde{p}_i^{\alpha}}\right] +\lambda_0\beta\widetilde{p}_i^{\beta-1}
+\beta\sum_{r=1}^m \lambda_r g_r(a_i)\widetilde{p}_i^{\beta-1}=0,~~~~~i=1, \ldots, n.
\label{EQ:MaxEnt1}
\end{eqnarray}
Multiplying (\ref{EQ:MaxEnt1}) by $\widetilde{p_i}$ and summing over all $i$, we get
$$
\lambda_0\beta +\beta\sum_{r=1}^m \lambda_r G_r=0,
~~\Rightarrow ~\lambda_0 = -\sum_{r=1}^m \lambda_r G_r,
$$
where we have used the conditions given in (\ref{EQ:Constraint_simple}).
Substituting the value of $\lambda_0$ in (\ref{EQ:MaxEnt1}), we get 
\begin{eqnarray}
\frac{\alpha\beta}{\alpha-\beta}\left[\frac{\widetilde{p}_i^{\beta-1}}{\sum_i \widetilde{p}_i^{\beta}}
- \frac{\widetilde{p}_i^{\alpha-1}}{\sum_i \widetilde{p}_i^{\alpha}}\right] 
+\beta\sum_{r=1}^m \lambda_r (g_r(a_i)-G_r)\widetilde{p}_i^{\beta-1}=0,~~~~~i=1, \ldots, n.
\label{EQ:MaxEnt2}\nonumber
\end{eqnarray}
Dividing by $\widetilde{p}_i^{\beta-1}$, using (\ref{EQ:Constraint_simple})
and using the transformation $\lambda_r \mapsto \alpha\lambda_r$, we get
\begin{eqnarray}
&&\frac{\widetilde{p}_i^{\alpha-\beta}}{\sum_i \widetilde{p}_i^{\alpha}}= 
1 + \left(\alpha - {\beta}\right)\sum_{r=1}^m \lambda_r (g_r(a_i)-G_r),
~~~~~i=1, \ldots, n.
\label{EQ:MaxEnt3}\nonumber\\
\Rightarrow &&
\widetilde{p}_i
\propto \left[1 + \left(\alpha - {\beta}\right)\sum_{r=1}^m \lambda_r (g_r(a_i)-G_r)\right]^{\frac{1}{\alpha-\beta}},
~~~~~i=1, \ldots, n.
\label{EQ:MaxEnt4}\nonumber
\end{eqnarray}
Therefore, the maximizer of the LNE over $P\in \Omega_n^\ast$ (with a fixed $n$)
subject to the normalized $q$-expectation constraints is also given by 
\begin{eqnarray}
{p}_i \propto \left[1 + \left(\alpha - {\beta}\right)\sum_{r=1}^m \lambda_r (g_r(a_i)-G_r)\right]^{\frac{1}{\alpha-\beta}},
~~~~~i=1, \ldots, n.
\label{EQ:MaxEnt_sp}\nonumber
\end{eqnarray}
By normalization, we get the required MaxEnt distribution over $P\in \Omega_n$ as given in the theorem. 

\noindent\textbf{Part (b):} \\
Next, let us consider the particular subclass of the LNE family at $\alpha=\beta$.
We can proceed as above to obtain the corresponding first order conditions 
from the objective function (\ref{EQ:Obj_Fun}) with $\alpha=\beta$
which simplifies to the form
\begin{eqnarray}
-\beta \frac{\ln\widetilde{p}_i}{\sum_i \widetilde{p}_i^{\beta}}
+\beta \frac{\sum_i\widetilde{p}_i^{\beta}\ln\widetilde{p}_i}{\left(\sum_i \widetilde{p}_i^{\beta}\right)^2} +\lambda_0
+\sum_{r=1}^m \lambda_r g_r(a_i)=0,~~~~~i=1, \ldots, n.
\label{EQ:MaxEnt1_0}
\end{eqnarray}
Multiplying this Equation (\ref{EQ:MaxEnt1_0}) by $\widetilde{p_i}^\beta$ and summing over all $i$, 
and using (\ref{EQ:Constraint_simple}), we again get
$$
\lambda_0 +\sum_{r=1}^m \lambda_r G_r=0,
~~\Rightarrow ~\lambda_0 = -\sum_{r=1}^m \lambda_r G_r.
$$
Hence, along with the constraints in (\ref{EQ:Constraint_simple}), 
we finally get 
\begin{eqnarray}
&&\ln\widetilde{p}_i - \sum_i\widetilde{p}_i^{\beta}\ln\widetilde{p}_i= \sum_{r=1}^m \lambda_r (g_r(a_i)-G_r),
~~~~~i=1, \ldots, n.
\label{EQ:MaxEnt3_0}\nonumber\\
\Rightarrow && \widetilde{p}_i \propto \exp\left[\sum_{r=1}^m \lambda_r (g_r(a_i)-G_r)\right],
~~~~~i=1, \ldots, n.
\label{EQ:MaxEnt4_0}\nonumber
\end{eqnarray}
Therefore, after proper normalization, the final MaxEnt distribution corresponding to the LNE subclass with $\alpha=\beta$
turns out to be the one given in (\ref{EQ:MaxEnt_final_0}) completing the proof.
\hfill{$\square$}

\bigskip
Note that the MaxEnt distribution in (\ref{EQ:MaxEnt_final}) corresponding to the LNE with $\alpha\neq \beta$ 
can also be expressed in terms of the deformed $q$-exponential function (Definition \ref{DEF:q-deformed}) as
\begin{eqnarray}
{p}_i^\beta &\propto & e_{\alpha/\beta} \left[\sum_{r=1}^m \theta_r (g_r(a_i)-G_r)\right],
~~~\mbox{ with }~\theta_r = -\beta\lambda_r\in\mathbb{R}, ~r=1, \ldots, m, ~i=1, \ldots, n.
\label{EQ:MaxEnt_fqExp}
\end{eqnarray}
Thus, the escort distribution of the maximizer of the LNE over a $\beta$-linear family belongs to the $(\alpha/\beta)$-power-law family
(Definition \ref{DEF:beta-power-law}).
This result clearly generalizes the corresponding maximum entropy results for the R\'{e}nyi entropy;
see Section \ref{SEC:renyi_world} for more precise correspondence between the MaxEnt  distributions of the R\'{e}nyi and LNE entropies 
with $\alpha\neq \beta$.


Further, by taking limit $\alpha\rightarrow\beta$ in the MaxEnt distribution in (\ref{EQ:MaxEnt_final}) or (\ref{EQ:MaxEnt_fqExp}),
we also obtain the MaxEnt distribution of the form (\ref{EQ:MaxEnt_final_0}) for the LNE with $\alpha=\beta$ (with a different parametrization).
These MaxEnt distributions corresponding to the LNE subfamily with $\alpha=\beta$  belong to the classical exponential family  of distributions,
related to the extensive statistics governed by the Shannon entropy. 
This interesting phenomenon can be justified  through the scale-invariance of the LNE family and the concept of escort distributions;
in particular, the underlying reason is that the LNE of a sub-probability distribution at $\alpha=\beta$ equals the Shannon entropy of 
its escort distribution and the non-extensive normalized constraints can be viewed as linear expectation constraints in the escort distributions.

\begin{remark}
The LNE subfamily with $\alpha=\beta$  is the first set of generalized (non-Shannon) entropies 
that produces the classical exponential-type MaxEnt distribution under non-extensive constraints (i.e, over a $\beta$-linear family).
\label{REM:LNE_spl}
\end{remark}

In summary, the proposed LNE family, apart from containing the first scale-invariant entropies, 
is also the first class of entropies producing both the classical exponential type 
as well as non-extensive $\beta$-power-law type  MaxEnt distributions
under the nonextensive framework, along with an additional tuning parameter in both the cases 
for modeling more complex structures of the physical or information systems.

\section{The Associated Divergence and the Underlying Geometry}

The cross entropy and the divergence (relative entropy) are the other important components of information theory 
when a prior guess is available for the distribution under study 
\cite{Sundaresan:2007,Sason/Verdu:2018,Sason:2018,Hanawal/Sundaresan:2010, Kosut/Sankar:2017, Lapidoth/Pfister:2019, Boztas:2014,Bunte/Lapidoth:2014,Kumar/Sundaresan:2015a,Kumar/Sundaresan:2015b}. 
In this section and the subsequent one, we will define and study the cross entropy and relative entropy measures 
associated with our proposed LNE family in (\ref{EQ:LN_entropy})--(\ref{EQ:LN_entropy0}).

\subsection{The LNCE and the LNRE Families}

Let us consider any two (sub-)probability distributions $P=(p_1, \ldots, p_n)^T$ and $Q=(q_1, \ldots, q_n)^T$ 
both belonging to $\Omega_n^\ast$ with fixed $W(P)=W(Q)=W$, say ($W=1$ for probability distributions).
From the relations between existing entropy measures and associated cross entropy measures,
a natural definition for the (asymmetric) cross entropy measure between $P,Q$ 
corresponding to the entropy measure LNE can be given by 
\begin{eqnarray}
\mathcal{CE}_{\alpha, \beta}^{LN}(P, Q) 
&:=& \frac{\alpha\beta}{\alpha-\beta}\left[\frac{1}{\alpha}\ln\left({\sum_i p_i^{\alpha} q_i^{\beta-\alpha}}\right)
-\ln||P||_\beta\right] ,
~~~ \alpha, \beta>0, \alpha\neq\beta.
\label{EQ:LN_Centropy}\\
\mathcal{CE}_{\beta, \beta}^{LN}(P,Q) 
&:=& \beta\frac{\sum_i p_i^\beta \ln(p_i/q_i)}{\sum_i p_i^\beta} - \beta\ln||P||_\beta,
~~~ \beta>0.
\label{EQ:LN_Centropy0}
\end{eqnarray}
Note that,  for the uniform prior $Q=U=(\frac{W}{n}, \ldots, \frac{W}{n})$, we have 
$$
\mathcal{CE}_{\alpha, \beta}^{LN}(P, U) = \beta \ln\left({n}/{W}\right)-\mathcal{E}_{\alpha, \beta}^{LN}(P), ~~~\mbox{ for all } \alpha, \beta>0.
$$
Hence a minimizer of the above cross entropies in (\ref{EQ:LN_Centropy})--(\ref{EQ:LN_Centropy0})
with respect to $P$, given a uniform prior (and any additional constraints),  
coincides exactly with the corresponding MaxEnt distribution of the LNE family (under the same set of constraints);
they are clearly the dual optimization problems of each other.
We will refer to these cross entropy measures (\ref{EQ:LN_Centropy})--(\ref{EQ:LN_Centropy0}) 
as the logarithmic $(\alpha, \beta)$-norm cross entropy or the LNCE in short.
Further, for the particular case of $\beta=1$ and $P\in\Omega_n$ (probability distribution),
the LNCE coincides with the R\'{e}nyi measure of directed-divergence \cite{renyi:1961}
having parameter $\alpha>0$, just as the LNE coincides with the R\'{e}nyi entropy for probabilities at $\beta=1$.
Thus, the LNCE provides a two-parameter generalization of the R\'{e}nyi divergence over $\Omega_n^\ast$.

Recall that the interesting relative $\alpha$-entropy \cite{Lutwak/etc:2005,Kumar/Sundaresan:2015a} defined in (\ref{EQ:RE_alpha})
is related to the R\'{e}nyi divergence (\ref{EQ:RE_renyi}) through  the escort distribution as 
$\mathcal{RE}_\alpha(P, Q) = \mathcal{D}_{\frac{1}{\alpha}}(P_\alpha, Q_\alpha)$ for any $P, Q\in\Omega_n$. 
Noting the similar relationship between our LNE  and the R\'{e}nyi entropy in (\ref{EQ:LNE_escort}), 
we may now define a new relative entropy (divergence) corresponding to the LNE via the relation
\begin{eqnarray}
\mathcal{RE}_{\alpha, \beta}^{LN}(P, Q) = \frac{1}{\beta}\mathcal{D}_{\frac{\beta}{\alpha}}(P_\alpha, Q_\alpha), ~~~~P, Q\in\Omega_n^\ast~,
\label{EQ:LNRE_renyi}
\end{eqnarray}
which produces a two-parameter generalization of the relative $\alpha$-entropy.  
Clearly $\mathcal{RE}_{\alpha, \beta}^{LN}(P, Q)$ defines a proper statistical divergence; 
its properties have been studied initially in \cite{Ghosh/Basu:2018}, where it was referred to it as the relative $(\alpha, \beta)$-entropy measure.
Here, to be consistent with our previous notations relating with the LNE, 
we will refer to it as the logarithmic $(\alpha, \beta)$-norm relative entropy  or the LNRE in short.
Clearly, the LNRE coincides with the relative $\alpha$-entropy of \cite{Kumar/Sundaresan:2015a} 
when $\beta=1$. For general $\beta>0$ it has the form 
\begin{eqnarray}
\mathcal{RE}_{\alpha, \beta}^{LN}(P, Q) &=& \frac{\alpha}{\beta(\beta-\alpha)} \ln\left[\sum_i \left(\frac{p_i}{||P||_\alpha}\right)^{\beta}
\left(\frac{q_i}{||Q||_\alpha}\right)^{\alpha-\beta}\right], ~~~~~~~~~~~\alpha\neq \beta.
\label{EQ:LNRE}\\
\mathcal{RE}_{\alpha, \alpha}(P, Q) &=& \lim\limits_{\beta\rightarrow \alpha} \mathcal{RE}_{\alpha, \beta}^{LN}(P, Q) 
= \sum_i  \ln(\frac{p_i}{q_i})\left(\frac{p_i}{||P||_\alpha}\right)^\alpha + \ln\left(\frac{||Q||_\alpha}{||P||_\alpha}\right). 
\label{EQ:LNRE_alpha=beta}
\end{eqnarray}
The term $\frac{1}{\beta}$ in (\ref{EQ:LNRE_renyi}) allows us to define the LNRE as a divergence measure for any real $\beta\neq 0$,
and then the limiting cases for $\beta\rightarrow 0$ also yield a divergence class as given by
\begin{eqnarray}
\mathcal{RE}_{\alpha, 0}^{LN}(P, Q) &=& \lim\limits_{\beta\rightarrow 0} \mathcal{RE}_{\alpha, \beta}^{LN}(P, Q) 
= \sum_i  \ln(\frac{q_i}{p_i})\left(\frac{q_i}{||Q||_\alpha}\right)^\alpha + \ln\left(\frac{||P||_\alpha}{||Q||_\alpha}\right), 
\label{EQ:LNRE_beta0}
\end{eqnarray}

An immediate interesting property of the LNCE and LNRE families over the general space $\Omega_n^\ast$
of sub-probabilities is their scale-invariance along the lines of the LNE. 
In particular, the LNRE is scale-invariant with respect to both of its arguments but the LNCE is so only in its  first 
argument; for any  $c_1, c_2,  \alpha, \beta >0$ and $P, Q\in\Omega_n^\ast$, we have   
\begin{equation}
\mathcal{RE}_{\alpha, \beta}^{LN}(c_1P, c_2 Q) =\mathcal{RE}_{\alpha, \beta}^{LN}(P, Q),
~~~~\mbox{ and }  ~~~ 
\mathcal{CE}_{\alpha, \beta}^{LN}(c_1 P, Q) =\mathcal{CE}_{\alpha, \beta}^{LN}(P, Q).   
\label{EQ:LNCE_RE_invariance}
\end{equation}
On the other hand, the LNCE and LNRE families are neither symmetric in the choice of $(\alpha, \beta)$, like the LNE, 
nor in their arguments. However, for any $\alpha, \beta>0$ and any $P, Q\in\Omega_n^\ast$, the LNRE satisfies 
\begin{equation}
\mathcal{RE}_{\alpha, \beta}^{LN}(P, Q) =\mathcal{RE}_{\alpha, \alpha-\beta}^{LN}(Q, P).   
\label{EQ:LNRE_Symm}
\end{equation}
Additionally, the LNCEs are (functionally) related to the LNREs via the relation
\begin{eqnarray}
\mathcal{RE}_{\alpha, \beta}^{LN}(P, Q) = \frac{1}{\beta} \mathcal{CE}_{\beta, \alpha}^{LN}(P, Q) +  \frac{\alpha}{\beta} \ln||Q||_\alpha,
~~~~~~~~~~\alpha, \beta>0~,  P, Q\in\Omega_n^\ast.
\label{EQ:LNCE_LNRE}
\end{eqnarray}

In view of the above relations and the positiveness of $\alpha, \beta$,  for a given prior distribution $Q$, 
the minimization of the LNCE $\mathcal{CE}_{\alpha, \beta}^{LN}(P, Q)$ with respect to $P$ is indeed 
the same as the minimization of the relative entropy $\mathcal{RE}_{\beta, \alpha}^{LN}(P, Q)$ in its first argument $P$;
but the same is not necessarily true  for their minimization in the second argument $Q$ for a given $P$.
These \textit{forward} and \textit{reverse projections} represent very important topics in information science \cite{Kumar/Sundaresan:2015a};
but before we derive them in Section \ref{SEC:projections}, 
let us first discuss the basic geometry of the LNCE and LNREs.


\subsection{Continuity and Convexity}

In view of the relation (\ref{EQ:LNCE_LNRE}), the geometric properties of the new LNCE measure as a function of $P$
are exactly the same as  those of $\mathcal{RE}_{\beta,\alpha}^{LN}(P, Q)$;  
the second has been studied in \cite{Ghosh/Basu:2018} under more general set-ups. 
We here summarize the continuity and convexity properties of these LNCE and LNRE measures with respect to both of their arguments.

\begin{theorem}[Continuity]
For any $\alpha, \beta>0$, the LNCE functional $\mathcal{CE}_{\alpha, \beta}^{LN}(P, Q)$ and the LNRE functional 
$\mathcal{RE}_{\alpha, \beta}^{LN}(P, Q)$ 
are both continuous in $P\in\Omega_n$ (when $Q$ is held fixed) and also in $Q\in\Omega_n$ (when $P$ is held fixed).
\label{THM:Cont_1}
\end{theorem}
\noindent\textbf{Proof: }
The continuity of the LNRE follows from Remarks 2 and 4 of \cite{Ghosh/Basu:2018}.
Then, the continuity of the LNCE follows from that of the LNRE and the $L_\alpha$ norm via the relation (\ref{EQ:LNCE_LNRE}).
\hfill{$\square$}\\

It is known \cite{VanErven/Harremos:2014} that the R\'{e}nyi divergence is always jointly quasi-convex in both of its arguments
and convex in its second argument; additionally it can be shown to be jointly convex when its order lies in $[0,1]$.
On the other hand, the relative $\alpha$-entropy in (\ref{EQ:RE_alpha}) is neither convex nor bi-convex; 
rather it is only known to be quasi-convex in its first argument \cite{Kumar/Sundaresan:2015a}. 
In the following theorem, we will present some general convexity results for the LNCE and LNRE 
as generalizations of the R\'{e}nyi divergence and the relative $\alpha$-entropy, respectively.

\begin{theorem}[Convexity]
For any given $\alpha, \beta>0$ and $Q\in\Omega_n$, 
the LNCE functional $\mathcal{CE}_{\alpha, \beta}^{LN}(\cdot, Q)$ and the LNRE functional $\mathcal{RE}_{\alpha, \beta}^{LN}(\cdot, Q)$,
as a function of their first arguments, are both quasi-convex over the set of escort distributions $\Omega_{n}^{(\beta)}$.
\label{THM:Convx_1}
\end{theorem}
\noindent\textbf{Proof: }
See \cite{Ghosh/Basu:2018} and the relation (\ref{EQ:LNCE_LNRE}).
\hfill{$\square$}

\begin{remark}
See the equivalence between the LNE and the corresponding LNCE/LNRE via Theorems \ref{THM:LNE_prop5} and \ref{THM:Convx_1}. 
The LNE is concave over the set of appropriate escort distributions and accordingly 
the associated LNRE/LNCE is quasi-convex also over the set of escort distributions. 
\end{remark}

\begin{theorem}
For any given $\alpha > \beta>0$ and $P\in\Omega_n$, 
the LNRE functional $\mathcal{RE}_{\alpha, \beta}^{LN}(P, \cdot)$,
as a function of it second argument, is quasi-convex over the set of escort distributions $\Omega_{n}^{(\alpha-\beta)}$.
\label{THM:Convx_2}
\end{theorem}
\noindent\textbf{Proof: }
See \cite{Ghosh/Basu:2018}.
\hfill{$\square$}

\bigskip
Note that, at the particular choice $\beta=\alpha-1$, the LNRE functional $\mathcal{RE}_{\alpha, \alpha-1}^{LN}$ 	
is quasi-convex in its second argument over $\Omega_n$. 
Particularly,  $\mathcal{RE}_{2, 1}^{LN}$ is quasi-convex in both of its arguments (keeping the other component fixed) over $\Omega_n$.

\begin{remark}
Note that, the LNRE is convex in its second argument (keeping the first argument fixed) over $\Omega_n^{(\alpha)}$ 
for the special case $\beta=1$, where it coincides with the R\'{e}nyi divergence of $\alpha$-escort distributions.
However, the convexity of the general LNRE in its second argument is still an open question; 
although we do not have a proof at present, 
we conjecture that they would also be convex over an appropriate set of probabilities  for all choices of $\beta>0$.
\end{remark}

\subsection{Generalized Pythagorean Property}

The Pythagorean theorem is an important result in information theory. 
Although the concept was originally initiated in ancient mathematics relating coplanar points via their $L_2$ distance,
the Pythagorean theorem with respect to the relative entropy (KLD) has profound applications across information science and statistics. 
For the KLD measure defined in (\ref{EQ:RE_KL}), a convex set $\mathbb{E}\subseteq \Omega_n$ and a given $Q\in\Omega_n\setminus \mathbb{E}$, 
it states that
$$
\mathcal{RE}(P, Q) \geq \mathcal{RE}(P, P^\ast) + \mathcal{RE}(P^*, Q),
$$ 
where $P^\ast = \arg\min\limits_{P\in \mathbb{E}}\mathcal{RE}(P, Q)$ is called the information projection (or, I-projection) of $Q$ onto $\mathbb{E}$. 
This version is also known as the information theoretic Pythagorean theorem.
Interestingly, however, the R\'{e}nyi divergence does not satisfy the above Pythagorean theorem;
but it has recently been shown that a R\'{e}nyi divergence of order $\alpha\in(0, \infty)$ satisfies it whenever $\mathbb{E}$ is $\alpha$-convex
\cite{VanErven/Harremos:2014}.
A generalization of the above Pythagorean relationship is also developed for the relative $\alpha$-entropy 
over convex sets of probability distributions \cite{Kumar/Sundaresan:2015a}.

Here we prove that our LNRE also has such a generalized  Pythagorean property over the $\alpha$-convex sets of probability distributions. 
For this purpose, given any $Q\in\Omega_n$ and $r>0$, let us define the lower level sets associated with the LNRE
as $B_{\alpha,\beta}(Q; r) =\left\{P\in\Omega_n : \mathcal{RE}_{\alpha, \beta}^{LN}(P, Q)< r \right\}$.
Clearly $B_{\alpha,\beta}(Q; r)$ is always a $\beta$-convex set by Theorem \ref{THM:Convx_1}.

\begin{theorem}[Pythagorean Property]
Fix $\alpha, \beta>0$  with $\beta\neq \alpha$ and $P^{(0)}, P^{(1)}, Q \in \Omega_n$.
Suppose that $P_{\beta, \lambda}$ denotes the $(\beta, \lambda)$-mixture of $P^{(0)}$ and $P^{(1)}$ for any $\lambda\in[0,1]$,  
 $\mathcal{L}_\beta\subseteq\Omega_n^{(\beta)}$ denotes the line segment joining their respective escort distributions 
$P_\beta^{(0)}$ and $P_\beta^{(1)}$ and $\mathcal{L} = \left\{ P \in\Omega_n : P_\beta\in\mathcal{L}_\beta \right\}$. 
Then, we have the following results. 
	\begin{enumerate}
		\item[a)]  Suppose $\mathcal{RE}_{\alpha, \beta}^{LN}(P^{(0)}, Q)$ and $\mathcal{RE}_{\alpha, \beta}^{LN}(P^{(1)}, Q)$ are finite.
		Then, $\mathcal{RE}_{\alpha, \beta}^{LN}(P_{\beta, \lambda}, Q) \geq \mathcal{RE}_{\alpha, \beta}^{LN}(P^{(0)}, Q)$ for all $\lambda\in[0,1]$,
		i.e., $\mathcal{L}$ does not intersect $B_{\alpha,\beta}(Q; \mathcal{RE}_{\alpha, \beta}^{LN}(P^{(0)}, Q))$ over $\Omega_n$,
		if and only if 
		\begin{eqnarray}
		\mathcal{RE}_{\alpha, \beta}^{LN}(P^{(1)}, Q) \geq \mathcal{RE}_{\alpha, \beta}^{LN}(P^{(1)}, P^{(0)}) + \mathcal{RE}_{\alpha, \beta}^{LN}(P^{(0)}, Q).
		\label{EQ:Pyth1}
		\end{eqnarray} 
		
		\item[b)]  Suppose $\mathcal{RE}_{\alpha, \beta}^{LN}(P_{\beta, \lambda}, Q)$ is finite for some fixed $\lambda\in(0,1)$.
		Then, $\mathcal{L}$ does not intersect $B_{\alpha,\beta}(Q; \mathcal{RE}_{\alpha, \beta}^{LN}(P_{\beta, \lambda}, Q))$
		if and only if 
		\begin{eqnarray}
		\mathcal{RE}_{\alpha, \beta}^{LN}(P^{(1)}, Q) 
		&=& \mathcal{RE}_{\alpha, \beta}^{LN}(P^{(1)}, P_{\beta, \lambda}) + \mathcal{RE}_{\alpha, \beta}^{LN}(P_{\beta, \lambda}, Q),
		\label{EQ:Pyth2.1}\\
		\mbox{and }~~\mathcal{RE}_{\alpha, \beta}^{LN}(P^{(0)}, Q) 
		&=& \mathcal{RE}_{\alpha, \beta}^{LN}(P^{(0)}, P_{\beta, \lambda}) + \mathcal{RE}_{\alpha, \beta}^{LN}(P_{\beta, \lambda}, Q).
		\label{EQ:Pyth2.2}
		\end{eqnarray} 
	\end{enumerate}
	\label{THM:Pythagorean}
\end{theorem}
\noindent\textbf{Proof: }
The proof follows by direct applications of Theorem 4 of \cite{Ghosh/Basu:2018}.
\hfill{$\square$}\\

Note that, at $\beta=1$, the above theorem coincides with Theorem 9 of \cite{Kumar/Sundaresan:2015a}
and provides its generalization for our LNRE with any $\beta>0$ utilizing the concept of escort distributions. 
The following corollary to this theorem proves that the Pythagorean inequality of the R\'{e}nyi divergences over a $\beta$-convex set
\cite[Theorem 14,][]{VanErven/Harremos:2014} exactly holds also for the case of our generalized class of LNRE measures.

\begin{corollary}[Pythagorean Inequality]
	\label{THM:Pythagorean2}
Fix $\alpha, \beta>0$  with $\beta\neq \alpha$,  $Q \in \Omega_n$ and a subset $\mathbb{E}\subseteq\Omega_n$.
Suppose that $\mathbb{E}$ is $\beta$-convex and $P^\ast = \arg\min\limits_{P\in\mathbb{E}}\mathcal{RE}_{\alpha, \beta}^{LN}(P, Q)$ exists.
Then, we have the Pythagorean inequality
\begin{eqnarray}
\mathcal{RE}_{\alpha, \beta}^{LN}(P, Q) \geq \mathcal{RE}_{\alpha, \beta}^{LN}(P, P^\ast) + \mathcal{RE}_{\alpha, \beta}^{LN}(P^\ast, Q)
~~~~\mbox{ for all }~~ P\in \mathbb{E}.
\label{EQ:Pyth_ineq}
\end{eqnarray} 
\end{corollary}
\noindent\textbf{Proof: }
Fix any $P\in \mathbb{E}$ and apply Theorem \ref{THM:Pythagorean} with  the given $Q$ and with $P^{(0)}=P^\ast$ and $P^{(1)}=P$.
Since $E$ is $\beta$-convex, for any $\lambda\in[0,1]$,  the $(\beta, \lambda)$-mixture  $P_{\beta, \lambda}$ of $P^{(0)}$ and $P^{(1)}$ 
belongs to $\mathbb{E}$, and hence it satisfies 
$\mathcal{RE}_{\alpha, \beta}^{LN}(P_{\beta, \lambda}, Q) \geq \mathcal{RE}_{\alpha, \beta}^{LN}(P^{(0)}, Q)$ 
by definition of $P^{(0)}=P^\ast$. Then, the desired result (\ref{EQ:Pyth_ineq}) follows from Part (a) of Theorem \ref{THM:Pythagorean}.
\hfill{$\square$}\\

\begin{remark}
Although the LNCE is related to the LNRE via the relation (\ref{EQ:LNCE_LNRE}),
a Pythagorean type property of LNCE is in general much difficult unless we restrict on probabilities having unit $L_\beta$ norm only.
In particular, the first part (i) of Theorem \ref{THM:Pythagorean} holds for the LNCE $\mathcal{CE}_{\alpha, \beta}$
(in place of $\mathcal{RE}_{\alpha, \beta}$) if $||P^{(0)}||_\beta = 1$ and its second part (ii) holds if $||P_{\beta, \lambda}||_\beta=1$.
Also, Corollary \ref{THM:Pythagorean2} holds for the LNCE if  $||P^\ast||_\beta = 1$ which trivially holds for $\beta=1$ (R\'{e}nyi divergence).
\end{remark}

\section{The Projection Rules: Minimizing the LNRE and LNCE}
\label{SEC:projections}

The minimization of the divergence measures in either of its arguments has prominent applications 
in information theory, statistics, physics, etc.  
They are often referred to as the projection rules when the minimization is done under some constrained set of probability distributions,
e.g., under Tsallis non-extensive constraints discussed in Section \ref{SEC:LNE_MaxEnt}. 
When the minimization is done with respect to the first argument, keeping the second argument fixed,
the resulting minimizer is called the \textit{forward projection}, which is formally defined for our LNRE divergences as follows. 

\begin{definition}[Forward Projection]\label{DEF:Forward_Proj}
Fix $\alpha, \beta>0$, $Q\in \Omega_n$  and consider a subset $\mathbb{E}\subset \Omega_n$
such that $\mathcal{RE}_{\alpha, \beta}^{LN}(P, Q)<\infty $ for at least one $P\in \mathbb{E}$.
Then, the forward projection of $Q$ onto $\mathbb{E}$ with respect to the LNRE is defined as 
\begin{eqnarray}
P^\ast = \arg\min\limits_{P\in \mathbb{E}} \mathcal{RE}_{\alpha, \beta}^{LN}(P, Q) \in \mathbb{E}.
\label{EQ:Forward_Proj}
\end{eqnarray}
\end{definition}

\begin{remark}
The forward projection with respect to the LNCE can be defined similarly. 
Importantly, the forward projections with respect to the LNRE and LNCE coincide 
via their relation in (\ref{EQ:LNCE_LNRE}).
\label{REM:LNRE_forward_LNCE}
\end{remark}

We can also minimize a divergence measure with respect to its second argument 
keeping the first argument fixed. It is then referred to as the \textit{reverse projection} 
which finds applications in robust statistical inference in the presence of outliers. 
It is well known that the principle of maximum likelihood estimation is an instance of reverse projection
with respect to the KLD measure, with the first argument being fixed at the empirical probability distribution obtained from sample observations.
For our general LNRE divergences, we formally define the reverse projection as follows. 

\begin{definition}[Reverse Projection]\label{DEF:Reverse_Proj}
Fix $\alpha, \beta>0$, $P\in \Omega_n$  and consider a subset $\mathbb{E}\subset \Omega_n$
such that $\mathcal{RE}_{\alpha, \beta}^{LN}(P, Q)<\infty $ for at least one $Q\in \mathbb{E}$.
Then, the reverse projection of $P$ onto $\mathbb{E}$ with respect to the LNRE is defined as 
\begin{eqnarray}
Q^\ast = \arg\min\limits_{Q\in \mathbb{E}} \mathcal{RE}_{\alpha, \beta}^{LN}(P, Q) \in \mathbb{E}.
\label{EQ:Reverse_Proj}
\end{eqnarray}
\end{definition}

\begin{remark}
The reverse projection with respect to the LNCE can be defined similarly. 
However, unlike the forward projection, the reverse projections with respect to the LNRE and LNCE can differ quite significantly.
Recalling that the LNRE is scale-invariant also in its second argument (and the LNCE is not)
and that the convexity result with respect to the second argument is available only for the LNRE as of now, 
we will only discuss the reverse projections with respect to the LNRE divergence (Theorem \ref{THM:reverse_exist}).
\end{remark}

\subsection{Existence and Uniqueness of the Projection Rules}

We now provide some sufficient conditions for the existence and uniqueness of the projection rules 
obtained based on our LNRE divergences. These are discussed in more detail by \cite{Ghosh/Basu:2018}
and so we will only present the results (without proofs) for the sake of completeness.

In particular, whenever the forward projection with respect to the LNRE exists, 
its uniqueness can be proved by the Pythagorean property in Theorem \ref{THM:Pythagorean}; 
this is described formally in the following theorem. The next theorem talks about the sufficient conditions for its existence.

\begin{theorem}[Uniqueness of Forward Projection]
Consider a set $\mathbb{E} \subset \Omega_n$ that is $\beta$-convex and fix $Q\in \Omega_n$. 
If a forward projection of $Q$ onto $\mathbb{E}$ with respect to the LNRE exists, it must be unique.
\label{THM:forward_unique}
\end{theorem}

\begin{theorem}[Existence of Forward Projection]
	Fix $\alpha, \beta>0$ with $\beta\neq \alpha$ and $Q\in \Omega_n$. 
	Given any set $\mathbb{E} \subset \Omega_n$, suppose that it is $\beta$-convex, closed with respect to the $L_{\alpha/\beta}$-norm
	and there exists $P\in \mathbb{E}$ satisfying $\mathcal{RE}_{\alpha, \beta}^{LN}(P, Q)< \infty$. 
	Then, a forward projection of $Q$ onto $\mathbb{E}$ with respect to the LNRE (or LNCE) always exists
	(which is also unique by Theorem \ref{THM:forward_unique}).	
	\label{THM:forward_exist}
\end{theorem}

We next present two results generalizing Theorem 10 of \cite{Kumar/Sundaresan:2015a} in case of our LNRE;
the first one presents a relation between the forward projection with respect to the LNRE 
and the generalized Pythagorean inequality in Corollary  \ref{THM:Pythagorean2}
and the second one proves the subspace-transitivity property of this forward projection onto $\beta$-convex sets.

\begin{theorem}\label{THM:Pythagorean_forward}
Fix $\alpha, \beta>0$ with $\beta\neq \alpha$ and $Q\in \Omega_n$. 
\begin{itemize}
	\item[a)] A probability distribution $P^\ast\in \mathbb{E} \cap B_{\alpha,\beta}(Q; \infty)$ 
is a forward projection of $Q$ onto a $\beta$-convex set $\mathbb{E}$
with respect to the LNRE, $\mathcal{RE}_{\alpha, \beta}$, if and only if 
every $P\in \mathbb{E} \cap B_{\alpha,\beta}(Q; \infty)$ satisfies the relation (\ref{EQ:Pyth_ineq}).
	\item[b)] If the $\beta$-escort distribution of the forward projection $P^\ast$ is an algebraic inner point of $\mathbb{E}^{(\beta)}$, 
then $\mathbb{E} \subset B_{\alpha,\beta}(Q; \infty)$ and equality holds in (\ref{EQ:Pyth_ineq}) for all $P\in \mathbb{E}$. 
\end{itemize}
\end{theorem}
\noindent\textbf{Proof: }\\
\textbf{Part(a):} The only if part is proved in Corollary \ref{THM:Pythagorean2}.
The if part holds trivially from the non-negativity of the LNRE measure. \\
\textbf{Part(b):} Fix any $P\in\mathbb{E}$. 
If the $\beta$-escort distribution of the forward projection $P^\ast$ is an algebraic inner point of $\mathbb{E}^{(\beta)}$, 
then $P^\ast$ is a $(\beta, \lambda)$-mixture of $P$ and $P^{(1)}$ for some $\lambda\in (0,1)$ and $P^{(1)}\in\mathbb{E}$
(by definition). Then, the equality in (\ref{EQ:Pyth_ineq}) holds by Part (b) of Theorem \ref{THM:Pythagorean} with $P=P^{(0)}$. 
\hfill{$\square$}

\begin{remark}
The above theorem can be used to easily prove Theorem \ref{THM:forward_unique}.
If $P_1^\ast$  and $P_2^\ast$ are two forward LNRE projection of some $Q\in\Omega_n$ onto a $\beta$-convex set $\mathbb{E}\subset\Omega_n$,
then $ \mathcal{RE}_{\alpha, \beta}^{LN}(P_1^\ast, Q) =  \mathcal{RE}_{\alpha, \beta}^{LN}(P_2^\ast, Q) $.
Also, by Part (a) of Theorem \ref{THM:Pythagorean_forward}, 
$ \mathcal{RE}_{\alpha, \beta}^{LN}(P_2^\ast, Q) \geq  \mathcal{RE}_{\alpha, \beta}^{LN}(P_2^\ast, P_1^\ast)  +  \mathcal{RE}_{\alpha, \beta}^{LN}(P_1^\ast, Q) $.
Therefore, we get  $\mathcal{RE}_{\alpha, \beta}^{LN}(P, Q) =0$ implying $P_1^\ast=P_2^\ast$ by the divergence property of the LNRE.
\end{remark}

\begin{theorem}[Subspace-Transitivity of the Forward LNRE Projection]~\\
Fix $\alpha, \beta>0$ with $\beta\neq \alpha$ and $Q\in \Omega_n$, 
and consider two $\beta$-convex subsets $\mathbb{E}_1$ and $\mathbb{E}$ of $\Omega_n$ such that $\mathbb{E}_1 \subset \mathbb{E}$. 
Let $P^\ast$ be a forward projection of $Q$ onto $\mathbb{E}$ with respect to the LNRE, $\mathcal{RE}_{\alpha, \beta}$,
and $P_1^\ast$ be the corresponding forward LNRE projection of $Q$ onto $\mathbb{E}_1$.
If the equality in (\ref{EQ:Pyth_ineq}) holds for every $P\in \mathbb{E}$, 
then $P_1^\ast$ is a forward LNRE projection of $P^\ast$ onto $\mathbb{E}_1$.
\label{THM:subspace-transitivity}
\end{theorem}
\noindent\textbf{Proof: }
Fix any $P\in\mathbb{E}_1\subset \mathbb{E}$. 
By part (a) Theorem \ref{THM:Pythagorean_forward} applied to $P_1^\ast$, 
we get 
\begin{eqnarray}
\mathcal{RE}_{\alpha, \beta}^{LN}(P, Q) \geq \mathcal{RE}_{\alpha, \beta}^{LN}(P, P_1^\ast) + \mathcal{RE}_{\alpha, \beta}^{LN}(P_1^\ast, Q).
\nonumber
\end{eqnarray} 
But, by the assumption that equality holds in (\ref{EQ:Pyth_ineq}), we get 
\begin{eqnarray}
\mathcal{RE}_{\alpha, \beta}^{LN}(P_1^\ast, Q) &=& \mathcal{RE}_{\alpha, \beta}^{LN}(P_1^\ast, P^\ast) + \mathcal{RE}_{\alpha, \beta}^{LN}(P^\ast, Q),
\nonumber\\ ~~~\mbox{ and }~~
\mathcal{RE}_{\alpha, \beta}^{LN}(P, Q) &=& \mathcal{RE}_{\alpha, \beta}^{LN}(P, P^\ast) + \mathcal{RE}_{\alpha, \beta}^{LN}(P^\ast, Q).
\nonumber
\end{eqnarray} 
Combining all three, we finally get 
\begin{eqnarray}
\mathcal{RE}_{\alpha, \beta}^{LN}(P, P^\ast) \geq \mathcal{RE}_{\alpha, \beta}^{LN}(P, P_1^\ast) + \mathcal{RE}_{\alpha, \beta}^{LN}(P_1^\ast, P^\ast).
\nonumber
\end{eqnarray} 
Since the above inequality holds for every $P\in\mathbb{E}_1$, again by Part (a) of Theorem \ref{THM:Pythagorean_forward}, 
it follows that $P_1^\ast$ is a forward LNRE projection of $P^\ast$ onto $\mathbb{E}_1$.
\hfill{$\square$}\\

Finally, in the next theorem, we present a set of sufficient conditions under which the reverse projection 
with respect to the LNRE exists and is unique; its proof follows directly from Theorem \ref{THM:forward_exist} 
via the relation (\ref{EQ:LNRE_Symm}).

\begin{theorem}[Existence and Uniqueness of Reverse Projection]
Fix $\alpha>\beta>0$ with $\beta\neq \alpha$ and $P\in \Omega_n$. 
Given any set $\mathbb{E} \subset \Omega_n$, assume that it is ${(\alpha-\beta)}$-convex,  closed 
and there exists $Q\in \mathbb{E}$ satisfying $\mathcal{RE}_{\alpha, \beta}^{LN}(P, Q)< \infty$. 
Then, a reverse projection of $P$ on $\mathbb{E}$ with respect the LNRE exists and is unique.
\label{THM:reverse_exist}
\end{theorem}


\subsection{Forward Projection onto a $\beta$-Linear Family}
\label{SEC:CE}

We now derive the constrained minimizer of the LNRE and the LNCE, in their first arguments, 
under the $\beta$-normalized (non-extensive) expectation constraints (\ref{EQ:Constraints_Nq}). 
This is clearly a dual problem to the constrained MaxEnt problem discussed in Section \ref{SEC:LNE_MaxEnt}. 
Further, for a given second argument $Q\in\Omega_n$, 
the solutions to the above constrained minimization problems are nothing but the forward projections 
obtained with respect to the LNRE  or the LNCE, respectively, for $Q$ onto a $\beta$-linear space defined as 
\begin{eqnarray}
\mathbb{L}_\beta = \left\{ P\in\Omega_n : \sum_{i} f_r(a_i)p_i^\beta = 0, ~~ r = 1, \ldots, m \right\}.
\label{EQ:beta-linear}
\end{eqnarray}
where $f_r = g_r - G_r$ for each $r=1, \ldots, m$. 
Here, we will derive the forward projections onto a general $\beta$-linear family generated by some $f_r$s; 
recall that these forward projections with respect to LNRE and LNCE are indeed the same (Remark \ref{REM:LNRE_forward_LNCE}).

Using the scale-invariance of the LNRE and the LNCE in their first argument, 
we can proceed as in the derivation of MaxEnt distribution in Section \ref{SEC:LNE_MaxEnt}
to derive the forward projections onto $\mathbb{L}_\beta$ via the method of Lagrange multiplier.
The resulting expressions of the projection distributions are given in the following theorem;
the proof is similar to that of Theorem \ref{THM:LNE_MaxEnt} and hence omitted.

\begin{theorem}
Fix $\alpha, \beta>0$, $Q=(q_1, \ldots, q_n)\in \Omega_n$ and consider a $\beta$-linear family of distributions $\mathbb{L}_\beta$
as defined in  (\ref{EQ:beta-linear}) for some given functions $f_1, \ldots, f_m$. 
Then the forward projection $P^\ast=(p_1^\ast, \ldots, p_n^\ast)\in\Omega_n$ of $Q$ onto $\mathbb{L}_\beta$,
obtained with respect to the LNRE (or the LNCE) having tuning parameters $(\alpha, \beta)$,   
has the following form. 
	\begin{itemize}
		\item[a)]  When $\alpha\neq \beta$, there exists constants $\lambda_1, \ldots, \lambda_m$ such that 
		\begin{eqnarray}
		{p}_i^\ast &=& \frac{1}{Z} \left[q_i^{\alpha-\beta} + \left(\alpha - {\beta}\right)\sum_{r=1}^m \lambda_r f_r(a_i)\right]^{\frac{1}{\alpha-\beta}},
		~~~~~i=1, \ldots, n,
		\label{EQ:LNRE_ForwardProj}\\
\mbox{with }~~~		Z &=& \sum_{i=1}^n \left[q_i^{\alpha-\beta} + \left(\alpha - {\beta}\right)\sum_{r=1}^m \lambda_r f_r(a_i)\right]^{\frac{1}{\alpha-\beta}}.
		\end{eqnarray}
		
		\item[b)] When $\alpha=\beta$, there exists constants $\lambda_1, \ldots, \lambda_m$ such that
		\begin{eqnarray}
		{p}_i^\ast &=& \frac{q_i\exp\left[\sum_{r=1}^m \lambda_r f_r(a_i)\right]}{
			\sum_i q_i\exp\left[\sum_{r=1}^m \lambda_r f_r(a_i)\right]},
		~~~~~i=1, \ldots, n.
		\label{EQ:LNRE_ForwardProj_0}
		\end{eqnarray}
In both cases, the $m$ constants  $\lambda_1, \ldots, \lambda_m$ are uniquely defined by 
the $m$ constraints defining the family $\mathbb{L}_\beta$.
	\end{itemize} 
	\label{THM:LNRE_ForwardProj}
\end{theorem}

Note that, the forward projection in the above theorem, obtained with respect to the LNRE (LNCE) measure with $\alpha\neq \beta$
onto a $\beta$-linear family has the same form as the forward projections of the KLD or R\'{e}nyi entropy onto a $\beta$-linear family \cite{Kumar/Sason:2016,Kumar/Sason:2016b}. It is also the same as the  the minimum Kullback-Leibler cross entropy
distribution obtained under nonextensive constraints \cite{Dukkipati/etc:2005}.
In fact, it can be easily seen that the forward projection in (\ref{EQ:LNRE_ForwardProj}) for the cases $\alpha\neq \beta$
belongs to an $(\alpha/\beta)$-power-law family generated by $Q_\beta$ and the functions $f_1, \ldots, f_m$.

\begin{remark}
Interestingly, the forward projection of the LNRE (LNCE) subclass at $\alpha=\beta$ belongs to the classical exponential family of distributions.
It is of the same form as the forward projection of KLD onto a linear family (i.e., under linear (expectation) constraint),
which is obtained here as a forward projection onto a $\beta$-linear family (under the non-extensive constraint).
As is consistent with Remark \ref{REM:LNE_spl} for the LNE, the subclass of the LNRE (or LNCE) family at $\alpha=\beta$ 
represents the first examples of divergence (or cross entropy) measures that have an exponential-type forward projection onto a $\beta$-linear space. 
\end{remark}

Let us now verify some natural axiomatic properties of the forward projection rule  of $Q\in\Omega_n$ onto a $\beta$-linear family $\mathbb{L}_\beta$, 
generated by the LNRE (or the LNCE) having tuning parameters $(\alpha, \beta)$, 
as a functional mapping ($\Pi_{\alpha, \beta}$) of $\mathbb{L}_\beta$ and $Q$ onto a value 
$\Pi_{\alpha, \beta}(\mathbb{L}_\beta|Q)\in \mathbb{L}_\beta$.
Clearly, $\Pi_{\alpha, \beta}(\mathbb{L}_\beta|Q)=P^\ast$ has the form as given in Theorem \ref{THM:LNRE_ForwardProj}.
We consider the axioms discussed by Csisz\'{a}r \cite{Csiszar:1991} 
while characterizing general projection rules and hence justifying inference procedures for linear families.   
However, here we will consider the $\beta$-linear families as mentioned above
and the associated results for the LNRE forward projection are given in the following theorem.

\begin{theorem}\label{THM:forwardProj_Lb}
Fix $\alpha, \beta>0$ with $\beta\neq \alpha$ and consider the corresponding forward projection rule (mapping) $\Pi_{\alpha, \beta}$
generated by the LNRE (or the LNCE) having tuning parameters $(\alpha, \beta)$.
\begin{itemize}
	\item[a)] The forward LNRE (or LNCE) projection rule   $\Pi_{\alpha, \beta}$ is \textit{regular} 
	\cite[][Definition 2]{Csiszar:1991}.
	\item[b)] If $\alpha<\beta$, then the forward LNRE (or LNCE) projection rule   $\Pi_{\alpha, \beta}$ is subspace-transitive
	\cite[][Definition 6]{Csiszar:1991} 	with respect to the $\beta$-linear families. 
\end{itemize}
\end{theorem}
\noindent\textbf{Proof: }
Part(a) follows from Theorem 1 of \cite{Csiszar:1991}.
Part(b) follows from Theorem \ref{THM:Pythagorean_forward} and \ref{THM:subspace-transitivity}
provided we can show that the $\beta$-escort distribution of the LNRE forward projection is an algebraic inner point of $\mathbb{L}^{(\beta)}$.
The last statement can be shown to hold for $\alpha<\beta$ 
by following the arguments similar to the proof of Proposition 15 in \cite{Kumar/Sundaresan:2015a},
with suitable modifications to replace their linear family by the $\beta$-linear family in the present case.  
\hfill{$\square$}

\begin{remark}
Csisz\'{a}r \cite{Csiszar:1991} proved that if any projection rule (with respect to the linear families) is 
regular, local \cite[][Definition 3]{Csiszar:1991} and subspace-transitive, 
then it must be generated by  (forward projection of) a Bregman divergences of sum-form \cite{Bregman:1967,Stummer/Vajda:2012}.
The KLD (\ref{EQ:RE_KL}) as well as the simple Euclidean distance falls under this category and satisfies all the three axioms. 
However, if we forgo the locality axiom, till date, there is only one known class of projection rules satisfying 
both the regularity and subspace-transitivity axioms with respect to the linear families (and not the locality axiom), 
which is generated by the relative $\alpha$-entropies of \cite{Kumar/Sundaresan:2015a}.
In Theorem \ref{THM:forwardProj_Lb}, we have now shown that the projection rules generated by our LNRE or LNCE 
also satisfy both the regularity and subspace-transitivity properties for the $\beta$-linear families.
This extends the findings of \cite{Kumar/Sundaresan:2015a} (they coincide when $\beta=1$) 
and motivates us to look  for such axiomatic characteristics of projection rules over linear as well as $\beta$-linear families. 
It would be an important future work to verify if these axioms are also satisfied 
by the forward projection rule generated by the LNRE onto linear families 
and if they can indeed characterize our LNRE family. 
\end{remark}

\section{Relation with R\'{e}nyi's Information World}
\label{SEC:renyi_world}

The R\'{e}nyi entropy and the associated R\'{e}nyi divergence have vast applications across several disciplines of (information) science. 
As the classical information theory is rooted to the Shannon entropy and the KLD, 
the world of modern information theory, generalized  based on the R\'{e}nyi entropy and the corresponding divergence measure,
has been observed to lead to significant improvements in several applications 
\cite{Yu/Tan:2018, Yu/Tan:2019, Sason/Verdu:2017, Sason/Verdu:2018, Bunte/Lapidoth:2014, Sason:2018, Hayashi/Tan:2016, Sason:2015, Nakiboglu:2018, Rached/etc:2001, Kuzuoka:2019}.

We have already noted in Section \ref{SEC:renyi_LNE} that the proposed LNE has a one-to-one correspondence with the R\'{e}nyi entropy 
via the $P\leftrightarrow P_\beta$ correspondence of the escort distributions. 
A similar correspondence also exists between our LNRE and the R\'{e}nyi divergence via the relation (\ref{EQ:LNRE_renyi})
which puts them in one-to-one correspondence again via $P\leftrightarrow P_\beta$ correspondence. 
We will now show that this correspondence indeed makes the underlying geometric properties, like Pythagorean property and projection rules, 
of our LNE and LNRE measures equivalent to those obtained for the R\'{e}nyi's information world.
These equivalences would further justify the potentiality of wider applications of our proposed information measures
in all the cases where the R\'{e}nyi's measure is used, with greater flexibility via the additional parameter present in our proposals.

\subsection{Comparisons of the MaxEnt Distributions}

We have proved in Section \ref{SEC:LNE_MaxEnt} that the maximum entropy distribution for the LNE with $\alpha\neq\beta$
under the Tsallis non-extensive constraints (equivalently under the $\beta$-linear family) belongs to 
the class of $(\alpha/\beta)$-power-law family of distributions. 
On the other hand it has already been proved \cite{Lenzi/etc:2000,Jizba/Arimitsu:2004a}  that the R\'{e}nyi entropy maximizer under the linear family 
(usual expectation constraints) is indeed an $\alpha$-power-law distribution.
So, there is a direct relation between the two MaxEnt distributions via the $P\leftrightarrow P_\beta$ correspondence
as presented in the following proposition.

\begin{proposition}
If $P$ maximizes the LNE with tuning parameter $(\alpha, \beta)$ under a $\beta$-linear family generated by functions $f_1, \ldots, f_m$,
and $\widetilde{P}$ maximizes the R\'{e}nyi entropy of order $(\alpha/\beta)$ under the linear family generated by the same  $f_1, \ldots, f_m$,
then there is an one-to-one relation between $P$ and $\widetilde{P}$ via the $P\leftrightarrow P_\beta$ correspondence.
In particular, $\widetilde{P}$ will be the $\beta$-escort distribution of $P$.
\end{proposition}
\noindent\textbf{Proof: }
The results follows from the relation between the LNE and the R\'{e}nyi entropy as noted in Section \ref{SEC:renyi_LNE}
and the one-to-one nature of the $P\leftrightarrow P_\beta$ correspondence.
\hfill{$\square$}\\

As a consequence of the above theorem, for any given $\beta$, 
the set of all MaxEnt distributions corresponding to the LNEs with varying $\alpha>0$ under a $\beta$-linear family of distributions
equals the set of all R\'{e}nyi entropy maximizers over the corresponding linear family of distributions
and hence is independent of the choice of $\beta$.

\subsection{Equivalence of the Forward Projections}

We now show that the forward projection problem for our LNRE is also equivalent to that for the R\'{e}nyi Divergence. 
For our LNRE, in Section \ref{SEC:projections}, we have considered the following forward projection problem:
\begin{itemize}
	\item Fix $\alpha, \beta>0$ with $\alpha\neq \beta$, $Q\in\Omega_n$ and $\mathbb{E}\subset \Omega_n$. Solve 
	\begin{eqnarray}
	\inf\limits_{P\in \mathbb{E}} \mathcal{RE}_{\alpha, \beta}^{LN}(P, Q). 
	\label{EQ:forward_problem_LNRE}
	\end{eqnarray} 
\end{itemize}
The existence and uniqueness of this projection problem have been investigated and we have seen that 
a set of sufficient condition for this purpose is that
\begin{eqnarray}
\mathbb{E} \mbox{ is $\beta$-convex  and closed in $L_{\alpha/\beta}$-norm}.
\label{EQ:forward_cond_LNRE}
\end{eqnarray}

On the other hand, for the R\'{e}nyi divergence, let us consider the corresponding  forward projection problem as follows:
\begin{itemize}
	\item Fix $\alpha >0$ with $\alpha\neq 1$, $Q\in\Omega_n$ and $\mathbb{E}_1\subset \Omega_n$. Solve 
	\begin{eqnarray}
	\inf\limits_{P\in \mathbb{E}_1} \mathcal{D}_{\alpha}(P, Q). 
	\label{EQ:forward_problem_renyi}
	\end{eqnarray} 
\end{itemize}
This problem has been solved in  \cite{Kumar/Sason:2016} in detail;
they developed a set of sufficient condition for the existence and uniqueness of this projection problem in (\ref{EQ:forward_problem_renyi})
which is 
\begin{eqnarray}
\mathbb{E}_1 \mbox{ is $\alpha$-convex  and closed in $L_1$-norm}.
\label{EQ:forward_cond_renyi}
\end{eqnarray}

The following theorem shows that the forward projection problems in (\ref{EQ:forward_problem_LNRE}) and (\ref{EQ:forward_problem_renyi})
are indeed equivalent via the $P\leftrightarrow P_\alpha$ correspondence.
 
\begin{theorem}
The forward projection problem (\ref{EQ:forward_problem_LNRE}) for the LNRE with a given $\alpha, \beta>0$, $\alpha\neq \beta$, 
is equivalent to the projection problem (\ref{EQ:forward_problem_renyi}) for the R\'{e}nyi entropy with $\alpha$ replaced by $(\beta/\alpha)$ 
and $\mathbb{E}_1$ replaced by $\mathbb{E}^{(\alpha/\beta)}$.
Further, the associated sets of sufficient conditions, namely the hypotheses in (\ref{EQ:forward_cond_LNRE}) and (\ref{EQ:forward_cond_renyi}),
are also identical under the $P\leftrightarrow P_\alpha$ correspondence.
\end{theorem}
\noindent\textbf{Proof: }
The first part follows from the relation (\ref{EQ:LNRE_renyi}) and the positiveness of $\beta$. 
The second part can be easily shown using Proposition \ref{PROP:escort}.
\hfill{$\square$}\\

\subsection{Equivalence of the Pythagorean Property}

We now describe the equivalence between the Pythagorean properties of our LNRE and the R\'{e}nyi divergences.
Let us fix $\alpha, \beta>0$  with $\beta\neq \alpha$,  $Q \in \Omega_n$ and a subset $\mathbb{E}\subseteq\Omega_n$.
Suppose that $\mathbb{E}$ is $\beta$-convex and the forward LNRE projection $P^\ast$ of $Q$ onto $\mathbb{E}$ exists.
Then, in Corollary \ref{THM:Pythagorean2},  we have shown that the LNRE measures satisfy the Pythagorean inequality (\ref{EQ:Pyth_ineq}).

Now, by Proposition \ref{PROP:escort}, one can easily show that if $\mathbb{E}$ is $\beta$-convex, 
then $\mathbb{E}^{(\alpha)}$ is $(\beta/\alpha)$-convex. Thus, using the relation (\ref{EQ:LNRE_renyi}),
we  get that $P_\alpha^\ast$ is the forward projection with respect to the R\'{e}nyi divergence of order $(\beta/\alpha)$ of $Q_\alpha$ 
onto the set $\mathbb{E}^{(\alpha)}$ and the projection satisfies the inequality
\begin{eqnarray}
\mathcal{D}_{\beta/\alpha}(P_\alpha, Q_\alpha) \geq \mathcal{D}_{\beta/\alpha}(P_\alpha, P_\alpha^\ast) 
+ \mathcal{D}_{\beta/\alpha}(P_\alpha^\ast, Q_\alpha) ~~~~\mbox{ for all }~~ P_\alpha\in \mathbb{E}^{(\alpha)}.
\label{EQ:Pyth_ineq_renyi}
\end{eqnarray} 
By virtue of the $P\leftrightarrow P_\alpha$ correspondence, the relation in (\ref{EQ:Pyth_ineq_renyi})
indeed recovers the Pythagorean inequality for the R\'{e}nyi divergence of order $(\beta/\alpha)$,
which was independently proved in \cite[][Theorem 14]{VanErven/Harremos:2014} and \cite[][proposition 1]{Kumar/Sason:2016}.

\section{Applications in Robust Statistical Inference}
\label{SEC:Applications}


Relative entropies, or, in statistical nomenclature, divergence measures, are very useful tools in statistical inference, 
as they have natural applications in parametric statistical estimation and testing of hypothesis. 
Given a statistical divergence (i.e., a relative entropy) and a parametric model, 
the technique of estimating the unknown parameter by selecting the model element ``closest'' to the data distribution 
(in terms of the given divergence) has obvious natural appeal and offers a bona-fide objective function 
which leads to the desired estimator through the relevant optimization procedure. 
For the relative entropy, where the first argument (distribution) is represented by the data distribution, 
the minimizer distribution indeed corresponds to the reverse projection of the relative entropy (divergence) measure.

Classical  maximum likelihood estimation procedure, the backbone of modern statistical theory,  
is also a minimum divergence method obtained  via the reverse projection of an appropriately estimated KLD measure.
As observed earlier, the correspondence between Fisher's likelihood theory of general statistical inference 
and Shannon's entropy and Shannon's information theory is now well established \cite{Kullback:1997}. 
The maximum likelihood estimator, the parameter estimate obtained by the minimization of the Kullback-Leibler divergence, 
has many theoretical optimality properties, and is the most efficient estimator under standard regularity assumptions. 
Modern theory of inference, however, recognizes that the conditions under which these optimal results are attained 
may not always be achieved in practice, 
and the violation of these assumptions (primarily through misspecification of the model and through the presence of noise) 
may lead to a drastic degradation in the performance of these otherwise efficient maximum likelihood estimators.

More generalized divergences have been used in statistical inference in recent times  to achieve greater stability (robustness) 
of the resulting inference under data contamination. 
In fact, statistical inference based on (different kinds of) divergences has grown to be a a lively research topic in present times, 
and, apart from the phenomenal growth of computing power, the primary reason for the same is the wide acceptance of the idea 
that not just efficiency, but the robustness of the procedure should also be a guiding factor in the selection of any inference procedure.  
Prominent research monographs (in the last three decades or so) which deal with statistical inference based on divergence measures 
include \cite{Pardo:2006,Vajda:1989,Liese/Vajda:1987,Basu/etc:2011}, 
all of which put some emphasis on the use of these divergences in the information theory context as well.

Our LNRE also has similar use in robust statistical inference 
leading to improved trade-offs between the efficiency and robustness of the resulting parameter estimates. 
Let $G$ be the true distribution under study, which has  a density $g$ with respect to some common dominating measure $\mu$.
We wish to model $G$ by a family of model (parametric) distributions 
$\mathcal{F}=\{ F_{\boldsymbol{\theta}} : \boldsymbol{\theta} \in \Theta \subseteq {\mathbb{R}}^p \}$, 
where the distributions within the set ${\cal F}$ are indexed by $\boldsymbol{\theta}$.
We assume that each $F_{\boldsymbol{\theta}}$ admits a density, say $f_{\boldsymbol{\theta}}$, with respect to $\mu$
which has the same support as $g$ (independent of $\boldsymbol{\theta}$).
We also assume that the distributional class ${\cal F}$ is identifiable, 
so that each  element $Q$ in ${\cal F}$ corresponds to an unique $\boldsymbol{\theta}\in  \Theta$. 
One way to estimate the \textit{best fitting functional} $\boldsymbol{\theta}$ is by minimizing 
a discrepancy measure quantified by the divergence or the relative entropy between $G$ and $Q$, where $Q$ varies over the elements of ${\cal F}$.  
Our LNRE fits in naturally in this application. 
The form of the LNRE, defined in (\ref{EQ:LNRE}) for discrete distributions, 
can easily be extended to the present set-up associated with a general dominating measure $\mu$ as in \cite{Ghosh/Basu:2018}; 
the LNRE between the true distribution $G$ and any element $Q$ in the model family $\mathcal{F}$ may be expressed as  
\begin{eqnarray}
\mathcal{RE}_{\alpha, \beta}^{LN}(G, Q) & = & \frac{\alpha}{\beta(\beta-\alpha)} \ln\left[\int \left(\frac{g}{||G||_\alpha}\right)^{\beta}
\left(\frac{q}{||Q||_\alpha}\right)^{\alpha-\beta}d\mu\right] 
\nonumber\\
&=& \frac{1}{\alpha-\beta} \ln\int g^{\alpha}d\mu 
- \frac{\alpha}{\beta(\alpha-\beta)} \ln\int g^{\beta}q^{\alpha-\beta}d\mu + \frac{1}{\beta}\ln\int q^{\alpha}d\mu, 
\label{EQ:Entropy_alpha_beta}
\end{eqnarray}
for any $\alpha\neq \beta$, where $q$ is the density function corresponding to $Q$. 
The LNRE for the cases $\alpha=\beta$ can be extended similarly as the limiting case.
Then, the \textit{best fitting parametric distribution} $P^\ast\in \mathcal{F}$ for the true distribution $G$, 
in the sense of minimum LNRE, is defined as the reverse projection of $G$ onto $\mathcal{F}$
assuming it exists and is unique.
Under the identifiability property of the model family, this reverse projection  $P^\ast\in \mathcal{F}$,
in turn, yields the functional ${\boldsymbol{\theta}}_{\alpha, \beta}(G)$ of the parameter $\boldsymbol{\theta}$ 
which satisfies $P^\ast = F_{{\boldsymbol{\theta}}_{\alpha, \beta}(G)}$. 
In other words, this minimum LNRE functional (MLNREF) of $\boldsymbol{\theta}$ can be defined as 
\begin{equation}
{\boldsymbol{\theta}}_{\alpha, \beta}(G) = \arg\min \limits_{\boldsymbol{\theta} \in \Theta } 
\mathcal{RE}_{\alpha, \beta}^{LN}(G,F_{\boldsymbol{\theta}}),
\label{EQ:MLNREE}
\end{equation}
whenever the minimum exists.
Sufficient conditions for the existence of the minimum divergence (minimum cross entropy) functional 
can be obtained from the corresponding results on reverse projection 
with respect to the LNRE (Theorem \ref{THM:reverse_exist}) along with the identifiability assumption of $\mathcal{F}$;
see Corollary 2 of \cite{Ghosh/Basu:2018}.

In practical applications involving real data, the true distribution $G$ is, of course, unknown, 
so a suitable substitute has to be used in Equation (\ref{EQ:MLNREE}). 
Although $G$ is unknown, we have the capability of choosing a random sample from the true distribution $G$. 
Using this random sample we can construct an empirical estimate $\widehat{G}$ of the true distribution  $G$ 
(or, rather, an estimate $\widehat{g}$ for the true density $g$)  which can then be plugged in place of $g$ in Equation (\ref{EQ:Entropy_alpha_beta}). 
As the associated distribution is discrete, one can estimate the density $g$
(which, in this case, is indeed a probability mass function) 
by the relative frequency of each point in the support computed on the basis of the given sample data.
Referring to this density as $\widehat{g}$, and representing the density $q$ as $f_{\boldsymbol{\theta}}$ 
(the density function corresponding to $F_{\boldsymbol{\theta}}$), 
the minimum LNRE estimator (MLNREE) is then obtained through the maximization of 
$$
\frac{\alpha}{\beta (\alpha - \beta)} \ln\int \widehat{g}^\beta 
f_{\boldsymbol{\theta}}^{\alpha - \beta}d\mu - \frac{1}{\beta} \int f_{\boldsymbol{\theta}}^\alpha d\mu,
$$
by ignoring the part of the objective function that has no role in the optimization 
and considering the equivalent maximization problem of the modified objective function. 

What explains the robust behaviour of the MLNREEs (or, at least, some members of the MLNREEs), 
which the maximum likelihood estimator fails to achieve? 
This has actually been indicated in the evolving literature for a while. In a statistical context, 
the two parameter relative entropy in Equation (\ref{EQ:Entropy_alpha_beta}) has been referred to as the LSD 
(or the logarithmic S-divergence) in \cite{Maji/etc:2014}, 
who constructed these divergences independently for the first time and used them for inference purposes. 
Asymptotic theoretical properties as well as the finite sample properties of these minimum divergence estimators have been studied 
in some detail in \cite{Maji/etc:2014,Maji/etc:2016,Ghosh/Basu:2018}. 
A particular one parameter sub-family of this minimum divergence  was proposed by \cite{Jones/etc:2001}, 
which has been referred to by many different names in the subsequent literature, 
but we will use the nomenclature logarithmic density power divergence (LDPD) in this connection.
The LDPD (between the distributions $G$ and $Q$) has the form 
\begin{eqnarray}
\frac{1}{\gamma} \ln \int g^{\gamma+1}d\mu  - \frac{\gamma +1}{\gamma} \ln \int g q^\gamma d\mu + \ln \int q^{\gamma +1} d\mu,
\label{EQ:LDPD}
\end{eqnarray}
and is obtained by replacing $\alpha = \gamma +1$ and $\beta = 1$ in Equation (\ref{EQ:Entropy_alpha_beta}). 
The relative $\alpha$-entropy of \cite{Kumar/Sundaresan:2015a}, defined in (\ref{EQ:RE_alpha}), 
is indeed also the same as this LDPD with $\alpha=\gamma+1$. 
The Kullback-Leibler divergence between $g$ and $f_{\boldsymbol{\theta}}$,
which is minimized by the maximum likelihood estimator under this set-up, 
is a member of the LDPD class for $\gamma \rightarrow 0$. 
The inference procedures resulting from this divergence generate estimating equations 
which provide a natural downweighting to anomalous observations by weighting the maximum likelihood score function by powers of densities. 
Replacing $q$ by $f_{\boldsymbol{\theta}}$ in Equation (\ref{EQ:LDPD}) and taking a derivative with respect to $\boldsymbol{\theta}$, 
we get the estimating equation 
\begin{eqnarray}
\frac{\int g f_{\boldsymbol{\theta}}^\gamma \boldsymbol{u}_{\boldsymbol{\theta}} d\mu}{\int g f_{\boldsymbol{\theta}}^\gamma  d\mu} 
= \frac{\int f_{\boldsymbol{\theta}}^{\gamma +1} \boldsymbol{u}_{\boldsymbol{\theta}} d\mu}{\int f_{\boldsymbol{\theta}}^{\gamma +1} d\mu}, 
\label{EQ:LDPD_EstEqn}
\end{eqnarray}
where $\boldsymbol{u}_{\boldsymbol{\theta}} = \frac{\partial}{\partial \boldsymbol{\theta}} \ln f_{\boldsymbol{\theta}}$,
the likelihood score function. 
The solution of Equation (\ref{EQ:LDPD_EstEqn}) corresponds to the minimum LDPD functional, 
whereas the solution of Equation (\ref{EQ:LDPD_EstEqn}) after replacing $g$ with $\widehat{g}$ leads to the minimum LDPD estimator. 
Notice that in Equation (\ref{EQ:LDPD_EstEqn}) the downweighting of the likelihood score function 
$\boldsymbol{u}_{\boldsymbol{\theta}}$ is done by the power of the density $f_{\boldsymbol{\theta}}^\gamma$. 
For the likelihood case, which corresponds to $\gamma = 0$, 
the downweighting  coefficient becomes identically equal to 1, so that there is no downweighting, 
explaining the failure of the likelihood based method in controlling anomalous observations. 
The form of Equation (\ref{EQ:LDPD_EstEqn}) indicates that the minimum LDPD estimator is an M-estimator, 
but it uses a model-dependent $\psi$ function, instead of the usual standard location-scale type $\psi$ functions. 

Notice that Equation (\ref{EQ:LDPD_EstEqn}) represents a normalized version of the estimating equation of 
the minimum density power divergence estimator of \cite{Basu/etc:1998}. 
It has been extensively argued in \cite{Fujisawa:2013} that estimators based on normalized equations may have many benefits 
which non-normalized estimating equations may not have. 
In particular such estimators appear to have small bias under heavy contamination compared to other estimators. 
Thus the minimum LDPD estimator enjoys advantages of having a normalized estimating equation, having density based downweighting 
and having a model based $\psi$ function in M-estimation.   

For the more general case of our LNRE, a similar exercise shows that the corresponding estimating equation is 
$$
\frac{\int g^\beta f_{\boldsymbol{\theta}}^{\alpha-\beta} \boldsymbol{u}_{\boldsymbol{\theta}} d\mu}{\int g^\beta f_{\boldsymbol{\theta}}^{\alpha-\beta}  d\mu }
= \frac{\int f_{\boldsymbol{\theta}}^{\alpha} \boldsymbol{u}_{\boldsymbol{\theta}} d\mu}{\int f_{\boldsymbol{\theta}}^{\alpha} d\mu}, 
$$
which is an obvious two parameter generalization of Equation (\ref{EQ:LDPD_EstEqn}). 
Notice that this estimating equation enjoys at least two of  the advantages indicated earlier, 
including downweighting by a power of the density (here $f_{\boldsymbol{\theta}}^{\alpha-\beta}$) and having a normalized estimating equation. 
Clearly the resulting MLNREE has strong robustness properties whenever $\alpha > \beta$.

As indicated, the statistical properties of these divergences are studied along with the theoretical investigations of the MLNREE 
(minimum LSD estimators) in \cite{Maji/etc:2014,Maji/etc:2016,Ghosh/Basu:2018} for discrete distributions. 
It has been observed that the asymptotic efficiency of the MLNREE becomes independent of $\beta$
and decreases slightly with increasing values of  $\alpha$.
The robustness properties of general MLNREEs have also been studied theoretically by higher-order influence functions and breakdown point analyses;
the robustness indeed increases as $\alpha$ increases or $\beta$ decreases.
The finite sample performance of the MLNREE has also been studied in \cite{Maji/etc:2014,Maji/etc:2016,Ghosh/Basu:2018}
through a variety of simulations and real data applications.
It has been observed that, for some particular choices of $(\alpha, \beta)$ depending on the situation, 
the new LNRE measure leads to parameter estimates that have better efficiency  and robustness trade-offs 
compared to  similar existing divergences or relative entropy measures. 

\begin{figure}[!h]
	\centering
	\subfloat[Box-plots under pure data from Bin(10, 0.1) with no outliers]{
		\includegraphics[width=\textwidth]{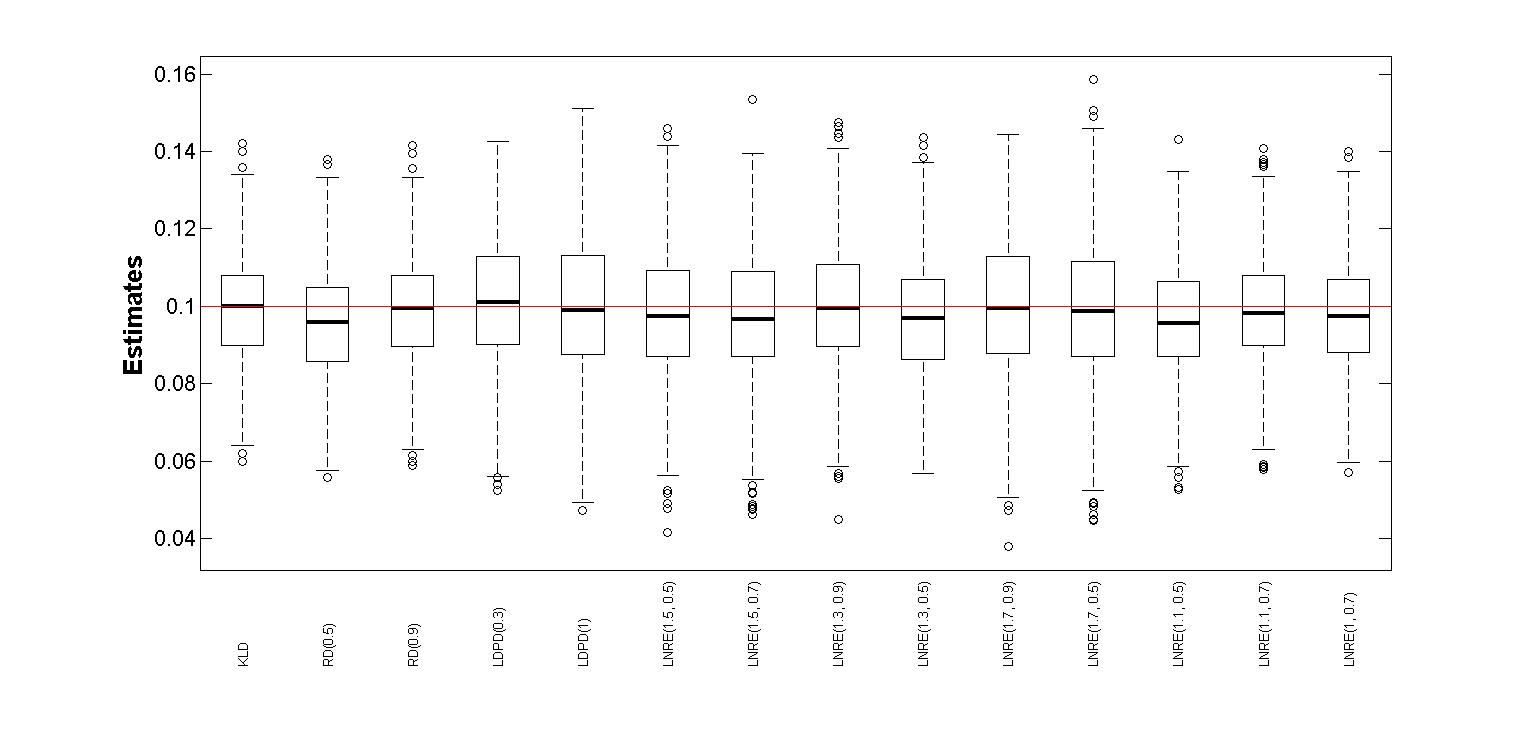}
		\label{FIG:pure_data}}
	\\
	\subfloat[Box-plots under contaminated data from Bin(10, 0.1) with 10\% outliers from Bin(10, 0.9)]{
		\includegraphics[width=\textwidth]{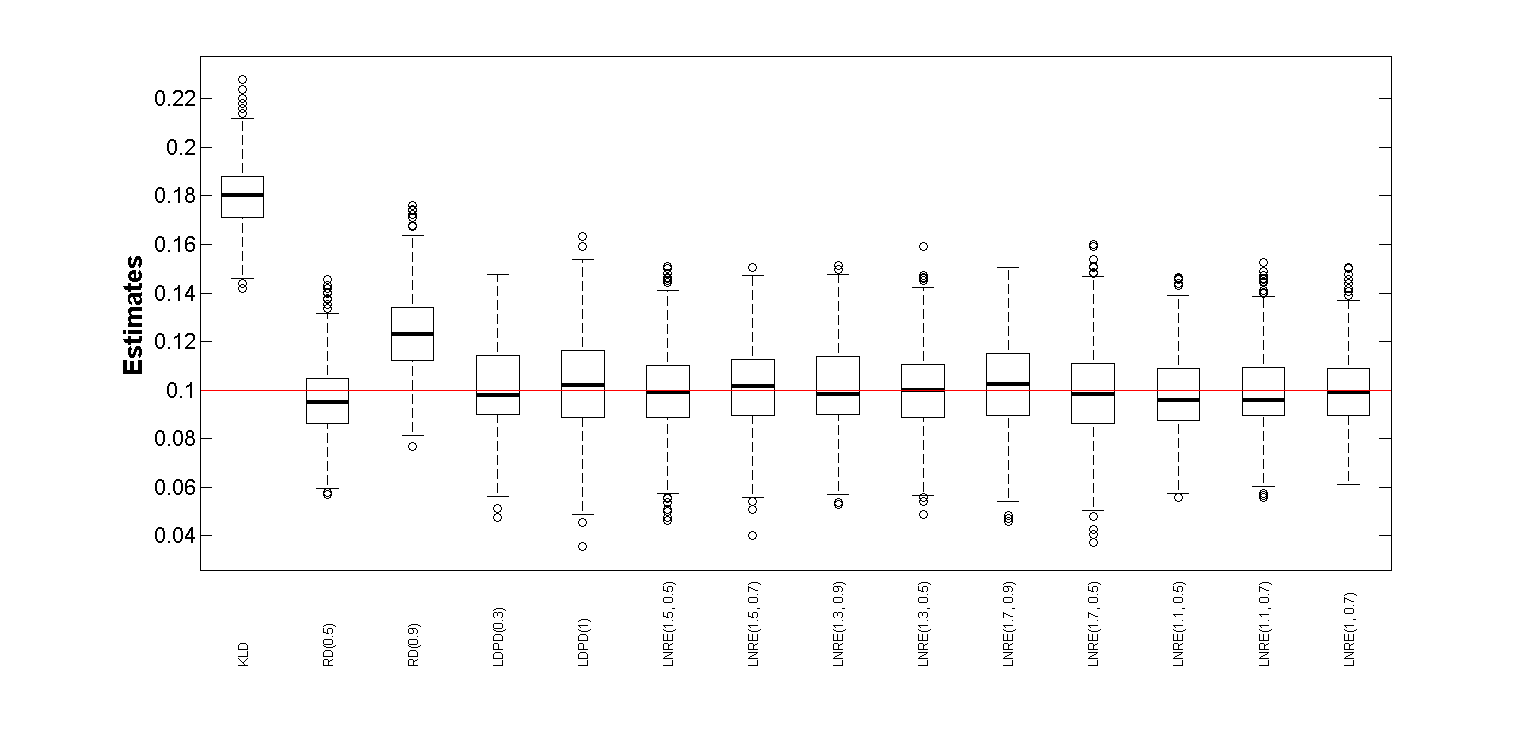}
		\label{FIG:contaminated data}}
	\caption{Box-plots of the parameter estimates obtained by minimizing different divergence measures 
		for pure and contaminated data from a binomial model. Here RD($\alpha$) denotes the R\'{e}nyi divergence of order $\alpha$,
	LDPD($\gamma$) denotes the LDPD with tuning parameter $\gamma$ and 
	LNRE($\alpha, \beta$) denotes our LNRE measure with tuning parameter $(\alpha, \beta)$}
	\label{FIG:boxplot}
\end{figure}

\bigskip\noindent
\textbf{Example 3. }[\textit{MLNREEs under the Binomial Model}]\\
As an illustrative example, let us compare the performance of the minimum LNRE estimators with some similar existing robust estimators
under the binomial model through a simulation study in the line of \cite{Ghosh/Basu:2018}.
Assuming that the true distribution $G$ is in $\Omega_{11}$ over the finite alphabet set $\mathcal{A}=\{a_0, a_1, \ldots, a_{10} \}$,
we can model it by the parametric family of binomial distributions, namely  $\mathcal{F}=\left\{ \mbox{Bin}(10, \theta) : \theta\in[0,1] \right\}$
and estimate the model parameter $\theta$ by minimizing appropriate divergence measures. 
We simulate a sample of size 50 from the true distribution  Bin(10, 0.1) and 
numerically compute the minimum LNRE estimator of $\theta$ based on this simulated sample for several choices of $\alpha>\beta>0$. 
We replicate the process 1000 times and the box-plots of the resulting MLNREEs are shown in Figure \ref{FIG:pure_data}. 
For comparison purposes, we also compute the classical MLE (obtained minimizing the KLD)
and the existing robust choices like the minimum R\'{e}nyi divergence estimators with $\alpha\in(0,1)$ and 
minimum LDPD (or relative $\alpha$-entropy) estimators with $\gamma\in(0,1)$ from the same set of samples;
the resulting estimates are also presented in Figure \ref{FIG:pure_data}. 
Next, to study the robustness properties of these minimum divergence estimators, 
we repeat the above simulation study by contaminating 10\% of each sample by independent observations simulated from Bin(10, 0.9) distributions.
The box-plots of the resulting estimates based on these contaminated samples are presented in Figure \ref{FIG:contaminated data}.
It can be clearly observed that, for some particular choices of $(\alpha, \beta)$ depending on the situation, 
the new LNRE measure leads to parameter estimates that have better efficiency  and robustness trade-offs 
compared to the similar existing divergences or relative entropy measures like the R\'{e}nyi divergence or the LDPD or the relative $\alpha$-entropy. 
\hfill{$\square$}
\\

%
The LNRE measure and the associated minimum LNRE estimates are also used in constructing robust testing procedures
for both one-sample and two sample hypothesis testing problems in \cite{Maji/etc:2014,Maji/etc:2016}.
Their properties have also been investigated in detail which reveals  similar advantages over the existing divergence based testing procedures,
further indicating the usefulness of our new LNRE measures.

\newpage
\section{Conclusions}

We believe that this paper has made several important contributions. 
We have proposed a new class of entropy measures and studied their properties in detail. 
This family has several  distinguishing features, the most remarkable of which is the scale-invariance property. 
This two parameter family of entropies contains the R\'{e}nyi entropy as a special case; 
the family can, in fact, be viewed as an appropriate R\'{e}nyi entropy of the escort distribution associated 
with a given sub-probability distribution. 
We have derived the MaxEnt distribution corresponding to our family of entropies  
under the Tsallis non-extensive constraints. 
This is in fact the first scale-invariant family of generalized entropies 
which provides the usual Maxwell-Boltzmann MaxEnt theorem of Shannon under nonextensive conditions. 
The relation and link between this research and our previous work in \cite{Ghosh/Basu:2018} have also been carefully explained. 
And finally, we have developed the corresponding cross entropies and 
studied some of their important properties including those of  their minimizers. 
We have also provided a sub-family of the cross entropy measures which leads to extensive results 
under the non-extensive framework, providing an instance of a case where the two concepts are combined. 

The present paper further opens up several interesting future research problems. 
The first one would be the application of our new entropy and cross entropy measures in generalizing
the concept of code-length or related information measures and the inference under source uncertainty or noisy channels.
It will also be interesting to generalize and study the statistical physics and thermodynamic concepts 
through the use of the new LNE family and its MaxEnt distributions; 
this will generalize the R\'{e}nyi thermodynamics and its applications.
Since it is already known that the minimizer of the associated  relative entropy measure leads to robust statistical inference \cite{Ghosh/Basu:2018},
it will also be interesting to study the applications of the LNE, its MaxEnt and the LNCE measures 
in different statistical and general inferential problems 
which are expected to provide robust solutions under data contaminations or outliers. 
We hope to pursue some of these extensions in our future work. 

\appendix

\section{Important Definitions}

\begin{definition}[$q$-deformed functions]
	\label{DEF:q-deformed}
	For any $q\in \mathbb{R}$, we define the $q$-\textit{logarithm} and $q$-\textit{exponential}
	of $x\in\mathbb{R}$, respectively, as
	\begin{eqnarray}
	\ln_q(x) &=& \frac{x^{1-q}-1}{1-q}, ~~~~ x>0, ~q\neq 1.
	\label{EQ:q-log}\\
	e_q(x) &=& \left\{\begin{array}{ll}
	[1+(1-q)x]^\frac{1}{1-q}, & \mbox{if } [1+(1-q)x]\geq 0,\\
	0 & \mbox{otherwise}
	\end{array}\right..
	\label{EQ:q-exp}
	\end{eqnarray}
	They coincide with the usual definitions of (natural) logarithm and exponential functions at $q\rightarrow 1$.
\end{definition}

\begin{definition}[$\beta$-Escort Distribution]
Given a (sub-)probability distribution $P\in \Omega_n^\ast$ and $\beta>0$, 
its escort probability distribution $P_\beta\in\Omega_n$ is defined as 
\begin{eqnarray}
P_\beta = \left(p_{1, \beta}, \ldots, p_{n, \beta}\right), ~~~~~\mbox{ with }~~
p_{i, \beta} = \frac{p_i^\beta}{||P||_\beta^\beta}, ~~i=1, \ldots, n.
\label{EQ:escort}
\end{eqnarray}
\label{DEF:escort}
\end{definition}

\begin{definition}[The $P \leftrightarrow P_\beta$ Correspondence]
Given a probability distribution $P\in\Omega_n$, 
a function $P_\beta$ is said to be in correspondence with $P$, 
denoted as $P \leftrightarrow P_\beta$, whenever they satisfy the relation (\ref{EQ:escort}).
\label{DEF:p-pbeta}
\end{definition}

\begin{definition}[$\beta$-linear family]
	For any fixed $\beta>0$ and $m$ given functions $f_1, \ldots, f_m$ on a finite alphabet set $\mathcal{A}=\{a_1, \ldots, a_n\}$, 
	the corresponding $\beta$-linear family of probability distributions on $\mathcal{A}$ is defined as
	\begin{eqnarray}
	\mathbb{L}_\beta = \left\{ P\in\Omega_n : \sum_{i} f_r(a_i)p_i^\beta = 0, ~~ r = 1, \ldots, m \right\}.
	\end{eqnarray}
	When $\beta=1$, the family $\mathbb{L} : = \mathbb{L}_1$ is referred to as the linear family. 
	\label{DEF:beta-linear}
\end{definition}

\begin{definition}[$\beta$-power-law family]
	Given any $\beta>0$, a probability distribution $Q\in\Omega_n$ (with full support when $\beta>1$), $m\in\{1, 2, \ldots \}$ and 
	$\Theta =\left\{ \theta=(\theta_1, \ldots, \theta_m) : \theta_r \in \mathbb{R}, r=1, \ldots, m \right\} \subset \mathbb{R}^m$,
	the $\beta$-power-law family generated by $Q$ and a set of $m$ given functions $g_r : \mathcal{A} \mapsto \mathbb{R}$, $r=1, \ldots, m$, 
	is defined as the set of probability distributions $P_\theta$, $\theta\in\Theta$, given by 
	\begin{eqnarray}
	P_\theta(x) &=& \frac{1}{Z_\beta(\theta)} \left[ Q(x)^{\beta-1} + (1-\beta) \sum_{r=1}^m \theta_r g_r(x)\right]^{\frac{1}{\beta-1}},
	~~~\beta\neq 1, ~x\in\mathcal{A},
	\label{EQ:beta-power-law}
	\\
	\mbox{or } ~~~~ P_\theta(x) ^{-1} 
	&=& {Z_\beta(\theta)} e_\beta\left[ \frac{Q(x)^{\beta-1}-1}{1- \beta} + \sum_{r=1}^m \theta_r g_r(x)\right], 
	~~~~~~~~\beta\neq 1, ~x\in\mathcal{A},
	\nonumber
	\end{eqnarray}
	where $Z_\beta(\theta)$ is the normalizing constant. 
	As $\beta\rightarrow 1$, this distribution coincides with the classical exponential family of distribution given by
	\begin{eqnarray}
	P_\theta(x) &=& \frac{1}{Z_1(\theta)} Q(x)\exp\left[ -\sum_{r=1}^m \theta_r g_r(x)\right], 
	~~~~x\in\mathcal{A}.
	\label{EQ:beta-power-law1}
	\end{eqnarray}
	\label{DEF:beta-power-law}
\end{definition}

\begin{definition}[$(\beta, \lambda)$-mixture]
Given two probability distributions $P_0, P_1 \in \Omega_n$, any $\beta>0$  and $0<\lambda <1$, 
the  $(\beta, \lambda)$-mixture of $P_0=(p_{01}, \ldots, p_{0n})$ and $P_1=(p_{11}, \ldots, p_{1n})$ is the probability distribution 
$P_{\beta, \lambda}\in\Omega_n$ defined as
\begin{eqnarray}
P_{\beta, \lambda}=(p_{\beta, \lambda}(i) : i=1, \ldots, n), ~~~~~\mbox{ with }
p_{\beta, \lambda}(i) =Z^{-1}\left[\lambda p_{0i}^\beta + (1-\lambda) p_{1i}^\beta \right]^{\frac{1}{\beta}}, ~~
i=1, \ldots, n,
\end{eqnarray}
where $Z$ is the required normalizing constant to make it into a probability distribution. 
It is important to note that $0<Z \leq 2$ for any non-trivial probability distribution $P_0, P_1$, as noted in \cite{Karthik/Sundaresan:2018},
and hence the above mixture distribution is well-defined.
\end{definition}

\begin{definition}[$\beta$-convex set]
	A subset $E\in\Omega_n$ is called $\beta$-convex for a given $\beta>0$ 
	if for any $P_0, P_1 \in \Omega_n$ and any $\lambda\in(0,1)$ their $(\beta, \lambda)$-mixture $P_{\beta, \lambda}$ belongs to $E$.
	\label{DEF:beta-convex}
\end{definition}

\begin{definition}[$\beta$-concave functional]
A functional $f$ from $\Omega_n$ to $\mathbb{R}$ is said to be $\beta$-concave for a given $\beta>0$ 
if for any $P_0, P_1 \in \Omega_n$ and any $\lambda\in(0,1)$ their $(\beta, \lambda)$-mixture $P_{\beta, \lambda}$ satisfies
\begin{eqnarray}
f(P_{\beta, \lambda}) \geq \lambda f((P_0)_\beta) + (1-\lambda) f((P_1)_\beta). 
\label{EQ:beta-concaveF}
\end{eqnarray}
The functional  $f$ is said to be $\beta$-convex if the inequality in (\ref{EQ:beta-concaveF}) is reversed. 
\label{DEF:beta-concaveF}
\end{definition}


\begin{definition}[Algebraic inner Point]
	Any probability distribution $P\in\mathbb{E}\subseteq\Omega_n$ is called an algebraic inner point of the set $\mathbb{E}$ 
	if for any $Q\in\mathbb{E}$, there exists $\widetilde{Q}\in \mathbb{E}$ and $\lambda\in(0,1)$ such that  
	$P= \lambda Q + (1-\lambda) \widetilde{Q}$.
\end{definition}

\begin{definition}[Generalized Additivity Property of an Entropy]
An entropy functional, say $\mathcal{E}(\cdot)$, defined over $\Omega_n^\ast$, 
the set of sub-probability distributions in (\ref{EQ:ProbSet_finite_gen}), is said to have the \textit{Generalized Additivity Property}
if there exists a strictly monotone and continuous function $g$ such that,
for any $P=(p_1, \ldots, p_n)\in \Omega_n^\ast$ and $Q=(q_1, \ldots, q_m)\in\Omega_m^\ast$
with $W(P) + W(Q) \leq 1$,  we have 
\begin{eqnarray}
\mathcal{E}(P\cup Q) 
= 	g^{-1}\left[\frac{W(P) g\left(\mathcal{E}(P)\right) + W(Q) g\left(\mathcal{E}(Q)\right)}{W(P)+W(Q)}\right], 
\label{EQ:GN_mean_additivity}
\end{eqnarray}
where $P\cup Q = (p_1, \ldots, p_n, q_1, \ldots, q_m)\in\Omega_{n+m}^\ast$.\\
For the particular choice $g(x)=ax+b$ for some real $a, b$, the relation (\ref{EQ:GN_mean_additivity}) simplifies to 
\begin{eqnarray}
\mathcal{E}(P\cup Q) 
= 	\frac{W(P) \mathcal{E}(P) + W(Q) \mathcal{E}(Q)}{W(P)+W(Q)}, 
\label{EQ:mean_additivity}
\end{eqnarray}
and it is referred to as the Additivity Property (or, Mean-Additivity property) of the entropy $\mathcal{E}(\cdot)$. 

The R\'{e}nyi entropy defined in (\ref{EQ:renyi_entropy_sp}) has this generalized additivity property (\ref{EQ:GN_mean_additivity})
for the choice $g(x) = 2^{(\alpha-1)x/c}$ with $c=\ln(2)$,
whereas the Shannon entropy defined in (\ref{EQ:Shannon_entropy_sp}) has the additivity property (\ref{EQ:mean_additivity}).
\label{DEF:Gen_additivity}
\end{definition}

\begin{definition}[Branching/Recursivity  Property of an Entropy]
	An entropy functional, say $\mathcal{E}(\cdot)$, defined over $\Omega_n$, 
	the set of probability distributions in (\ref{EQ:ProbSet_finite}), is said to satisfy the \textit{Branching or the Recursivity  Property}
	if, for any $P=(p_1, \ldots, p_n)\in \Omega_n$ and $Q=(q_1, \ldots, q_m)\in\Omega_m$,  we have 
$$
\mathcal{E}(p_1, p_2, \ldots, p_{n-1}, p_nq_1, p_nq_2, \ldots, p_n q_m) 
= \mathcal{E}(p_1, \ldots, p_n) + p_n\mathcal{E}(q_1, \ldots, q_m).	
$$
It is satisfied by the Shannon entropy defined in (\ref{EQ:Shannon_entropy}).
\label{DEF:Branching}
\end{definition}


\end{document}